\documentclass[11pt]{amsart}
\usepackage{amsmath,amsthm, amscd, amssymb, amsfonts,color,mathrsfs}
\usepackage{hhline}
\usepackage[active]{srcltx}
\usepackage[all]{xy}
\usepackage{enumitem}
\usepackage{graphicx,color}

\newcommand{\letra}{\mathbf{y}}
\newcommand{\letrag}{\mathbf{g}}
\newcommand{\Gletra}{\mathbf{G}}
\newcommand{\ydG}{{}^{\ku \Gletra}_{\ku \Gletra}\mathcal{YD}}

\newcommand{\spec}{\operatorname{spec}}
\newcommand{\GK}{\operatorname{GKdim}}

\newcommand{\cB}{\mathcal{B}}
\newcommand{\W}{{\mathcal W}}

\newcommand{\cC}{\mathcal{C}}
\newcommand{\cD}{\mathcal{D}}
\newcommand{\Ec}{\mathcal{E}}
\newcommand{\Fc}{\mathcal{F}}
\newcommand{\cG}{\mathcal{G}}
\newcommand{\cH}{\mathcal{H}}

\newcommand{\ku}{\Bbbk}
\newcommand{\I}{{\mathbb I}}
\newcommand{\N}{{\mathbb N}}
\newcommand{\G}{{\mathbb G}}

\newcommand{\U}{{\mathcal U}}
\newcommand{\Oc}{{\mathcal O}}

\newcommand{\Vc}{{\mathcal V}}
\newcommand{\Ac}{{\mathcal A}}

\newcommand{\ydg}{{}^{\ku G}_{\ku G}\mathcal{YD}}
\newcommand{\ydH}{{}^{H}_{H}\mathcal{YD}}
\newcommand{\Rc}{{\mathcal R}}

\newcommand{\jordan}{{\mathtt{J}}}

\newcommand{\hopfa}{\textsf{A}}
\newcommand{\hopfb}{\textsf{B}}

\newcommand{\End}{\operatorname{End}}

\newcommand{\Gln}{\operatorname{GL}}

\newcommand{\rep}{\operatorname{rep}}
\newcommand{\reps}{\operatorname{rep}_{\textrm{ssi}}}

\newcommand{\nuc}{\operatorname{ker}}
\newcommand{\dete}{\operatorname{det}}

\newcommand{\la}{\langle}
\newcommand{\ra}{\rangle}

\newcommand{\leftaction}{\mathfrak L}
\newcommand{\rightaction}{\mathfrak R}

\newcommand{\Kg}{\mathscr K}
\newcommand{\Igo}{\mathscr I}
\newcommand{\Hg}{\mathscr H}

\newcommand{\toba}{\mathscr B}
\newcommand{\base}{{\mathfrak{B}}}

\numberwithin{equation}{section}\theoremstyle{plain}

\newtheorem{theorem}{Theorem}[section]
\newtheorem{lema}[theorem]{Lemma}
\newtheorem{cor}[theorem]{Corollary}

\newtheorem{prop}[theorem]{Proposition}

\theoremstyle{definition}

\theoremstyle{remark}
\newtheorem{obs}[theorem]{Remark}

\newtheorem{step}{Step}
\newtheorem{passo}{Step}

\newcommand{\id}{\operatorname{id}}
\newcommand{\Id}{\operatorname{Id}}
\newcommand{\Ig}{\operatorname{Im}}
\newcommand{\ext}{\operatorname{Ext}}

\newcommand{\GL}{\mathbf{GL}}

\def\pf{\begin{proof}}
\def\epf{\end{proof}}

\theoremstyle{remark}

\begin{document}

\thispagestyle{empty}

\title[On the Bosonization of the Super Jordan  Plane]{On the Bosonization of the Super Jordan  Plane}
\author[Andruskiewitsch, Bagio, Della Flora, Fl\^ores]
{Nicol\'as Andruskiewitsch, Dirceu Bagio, Saradia Della Flora, Daiana Fl\^ores}

\address{FaMAF-Universidad Nacional de C\'ordoba, CIEM (CONICET),
Medina A\-llen\-de s/n, Ciudad Universitaria, 5000 C\' ordoba, Rep\'ublica Argentina.} \email{andrus@famaf.unc.edu.ar}

\address{Departamento de Matem\'atica, Universidade Federal de Santa Maria,
97105-900, Santa Maria, RS, Brazil} \email{bagio@smail.ufsm.br, saradia.flora@ufsm.br, flores@ufsm.br}

\thanks{
N. A. was partially supported by CONICET, MATHAMSUD, Secyt (UNC).\\
D. B. and D. F. were supported by FAPERGS 2193-25.51/13-3, MATHAMSUD}

\dedicatory{Dedicated to Professor Ivan Shestakov on the occasion of his 70th birthday}

\begin{abstract}
Let $H$ and $K$ be the bosonizations of the Jordan and super Jordan plane by the group algebra of a cyclic group; the algebra $K$ projects onto an algebra $L$ that can be thought of as the quantum Borel of $\mathfrak{sl}(2)$ at $-1$. The finite-dimensional simple modules over $H$ and $K$, are classified; they all have dimension 1, respectively $\le 2$.
The indecomposable $L$-modules of dimension $\leq 5$ are also listed.
An interesting monoidal subcategory of $\rep L$ is described.
\end{abstract}

\maketitle

\setcounter{tocdepth}{1}

\tableofcontents

\section{Introduction}

In this paper we start the study of the representation theory of two Hopf algebras introduced in
\cite{aah-jordan} and \cite{aah-triang}. In this last paper, Nichols algebras
over abelian groups with finite Gelfand-Kirillov dimension  were classified under some suitable hypothesis.
The list includes the well-known Jordan plane denoted here $\Ac$ 
and a new Nichols algebra called the super Jordan plane denoted here $\cB$. 
The study of various aspects of $\cB$ was undertaken in \cite{rep-super,RS}.
Let $G$ be an infinite cyclic group denoted multiplicatively with a fixed generator $g$.
The Nichols algebras $\Ac$ and $\cB$ are realized in $\ydg$ and we have then the Hopf algebras
$H = \Ac \# \ku G$ and $K = \cB \# \ku G$, where $\#$ stands for the Radford-Majid bosonization.
The liftings or deformations of $H$ and $K$ were computed in \cite{aah-jordan}.

The purpose of this paper is to begin the study of the tensor category of $K$-modules.
For this it is useful to study the $H$-modules, since $H$ embeds into $K$, and also the $L$-modules,
where $L$ is the quantum Borel of $\mathfrak{sl}(2)$ at $-1$ that appears as a quotient of $K$.
We obtain:
\begin{itemize}[leftmargin=*]
\item The  classification of the simple objects in  $\rep K$, that reduces to those in $\rep L$,
see Proposition \ref{prop:simple_boso_super} and Theorem \ref{theor:simple_boso_super};
all have dimension 1 or 2.

\item The  classification of the simple objects in  $\rep H$, all of dimension 1,  and of those indecomposable of dimension 2, see Propositions \ref{prop:simple-H} and \ref{prop:pol_min}.

\item The  classification of the indecomposable objects in  $\rep L$ of dimensions $\leq 5$, see Propositions \ref{prop:L-dim2}, \ref{prop:L-dim3}, \ref{prop:L-dim4} and \ref{prop:L-dim5}. 

\item An interesting subcategory of $\rep L$, see Subsection \ref{subsec:tensor-simples-L}.

\end{itemize}

\subsection{Notations and conventions} \label{subsec:notation}
We denote the natural numbers by $\N$ and $\N_0=\N\cup 0$.
If $k < t \in \N_0$, then we denote $\I_{k, t} = \{n\in \N_0: k\le n \le t \}$,
and   $\I_{t} := \{1,\dots, t\}$.

We work over an algebraically closed field $\ku$ of characteristic 0.
The group of $n$-th roots of 1 in $\ku$ is denoted $\G_n$; $\G'_n$ is the subset of the primitive ones.

We denote by $e_1, \dots,  e_{n}$ the canonical basis of $\ku^n$, which is identified 
with the space of column vectors.
Given a vector space $V$,   $T\in \End V$ and  
 $\lambda \in \ku$, we denote $V_T^{\lambda} = \ker (T - \lambda)$ and $V_T^{(\lambda)} = \cup_{j \in \N} \ker (T - \lambda)^j$.
The Jordan block of size $r$ associates to  $\lambda$ is denoted by $\jordan_{r}(\lambda)$.
The set of eigenvalues of $T$  is denoted $\spec T$.

We write $X \leq Y$ to express that $X$ is a sub-object of $Y$ in a category $\mathfrak C$.
All modules are left modules. As usual, $\rep A$ is the category of finite-dimensional representations of an algebra $A$; we use indistinctly the languages of representation and module theories.  
The braided tensor category of left Yetter-Drinfeld modules over a Hopf algebra $H$ is denoted by $\ydH$.

\section{The  quantum Borel subalgebra of $\mathfrak{sl}(2)$ at $-1$}\label{sec:borel}
Let $\Rc=\ku[y]$ be the polynomial algebra in one variable.
We realize $\Rc$ as a Hopf algebra in $\ydg$ by declaring $\Rc_{g^n} = \ku y^n$, $g \cdot y^n = (-1)^ny^n$ and $y$ to be primitive. 
Let $L= \Rc \#\ku G$ the bosonization of $\Rc$ by $\ku G$; i.~e. $L$ 
is the algebra generated by $y$ and $g^{\pm 1}$ with relations $g^{\pm 1}g^{\mp 1} = 1$,
\begin{align}\label{rel:superbosox1}
gyg^{-1}=-y.&
\end{align}
The comultiplication, the counit and the antipode of $g$ and $y$ are determined  by
\begin{align*}
\Delta(g) &= g\otimes g,& \Delta(y) &= y\otimes 1+g\otimes y.
\end{align*}
The Hopf algebra $L$ is the quantum Borel subalgebra of $\mathfrak{sl}(2)$ at $-1$.
Let $L^{(n)} = L / \langle g^{2n} - 1\rangle$, a quotient Hopf algebra of $L$, $n \in \N$.
The subalgebra $\ku \langle g^2, y^2 \rangle$ is a central Hopf subalgebra, 
the ideal $I$ generated by $g^2 - 1$ and $y^2$ is a Hopf ideal and $L/I$ is the 4-dimensional Sweedler algebra $H_4$.
Thus we have  exact sequences of Hopf algebras
\begin{align*}
&\ku \langle g^2, y^2 \rangle \hookrightarrow L \twoheadrightarrow H_4,&
&\ku \langle g^{2n} \rangle \hookrightarrow L \twoheadrightarrow L^{(n)}.
\end{align*}
The Hopf algebras $L$ and $L^{(n)}$ are pivotal and $L^{(1)}$ is spherical.

\subsection{Simple $L$-modules}\label{subsec:simple1}
We state here the (probably well-known) classification of the finite-dimensional simple $L$-modules;
we give a proof for completeness.
Throughout this subsection, $V\in \rep L$.

Let  $a\in \ku^{\times}$ be an eigenvalue of $g$ and $V^{a} = V_{g}^{a}$.
Then $yV^{a}\subset V^{-a}$, thus $V^{a}\oplus V^{-a}$ is a nonzero $L$-submodule of $V$.
Also, $y V_{g}^{(\pm a)}\subset V_{g}^{(\mp a)}$.

We describe the one-dimensional representations of $L$.
Given $a\in \ku^{\times}$, let $\ku_a = \ku$  with the action
$y \cdot 1 = 0$, $g \cdot 1=a$.
Clearly $\ku_a \in \rep L$, $\ku_a\simeq \ku_b$ if and only if $a=b$,
and every one-dimensional representation of $L$ is like this.

Next we describe the two-dimensional irreducible representations of $L$.
Given $a,b\in \ku^{\times}$, let $\U_{a,b} = \ku^2$ with the representation
\begin{align} \label{simple:2dim}
g &\mapsto \begin{pmatrix}
a & 0 \\
0 & -a \end{pmatrix},&  y&\mapsto \begin{pmatrix}
0 & b \\
1 & 0 \end{pmatrix}.
\end{align}
It is easy to check that $\U_{a,b}\in \rep L$ is irreducible.

\begin{prop}\label{prop:simple_boso_super}
Let $V$ be a finite-dimensional simple $L$-module. Then either $V\simeq \ku_a$ for a unique $a\in \ku^{\times}$, or else  $V\simeq \U_{a,b}$ for some $a,b \in \ku^{\times}$. Moreover, $\U_{a,b}\simeq \U_{a,c}$ if and only if $b=c$.
\end{prop}

\pf First we prove that $\dim V\leq 2$.
Let $a$ be an eigenvalue of $g$. Since $V^{a}$ is invariant by $g$ and $y^2$ and  $gy^2=y^2g$, there exists $v \in V^{a}$ which is a common eigenvector of $g$ and $y^2$. Then the subspace $\langle v,yv\rangle$ of $V$ is a non-trivial $L$-submodule of $V$. Thus $V = \langle v,yv \rangle$. 

Assume that $\dim V = 2$. We claim that  $y$ is invertible.
Indeed, $\ker y$ is an $L$-submodule of $V$. If $\ker y=V$
and $v$ is an eigenvector  of $g$, then $\langle v\rangle \leq V$, a contradiction. Thus $\ker y=0$.
Now $g$ has exactly two eigenvalues $a$ and $-a$.
If $v\in V^{a}$, $v\neq 0$, then $yv\in V^{-a}$ and $\{v,yv\}$ is a basis of $V$. Suppose that $y^2v=bv+b'yv$, where $b,b'\in \ku$. By \eqref{rel:superbosox1} $b'=0$. Hence  $V\simeq \U_{a,b}$. It is straightforward to check that $\U_{a,b}\simeq \U_{a,c}$ if and only if $b=c$.
\epf

\begin{obs}
As is well-known, Proposition \ref{prop:simple_boso_super} extends to $q \in \G'_n$ instead of $-1$. Indeed,
let $L_q = \ku \langle y,g^{\pm 1} \vert g^{\pm 1}g^{\mp 1} - 1, gyg^{-1} - qy\rangle$, 
the quantum Borel of $\mathfrak{sl}(2)$ at $q$. Then any  simple object of $\rep L_q$ 
is isomorphic either to $\ku_a$, or else to  $\U_{a,b}$ for some $a,b \in \ku^{\times}$. Here 
$\ku_a = \ku$  with $y \cdot 1 = 0$, $g \cdot 1=a$; and 
$\U_{a,b} = \ku^n$ with $g\cdot e_i = aq^{i-1} e_i$, $i\in \I_n$, 
$y\cdot e_j =  e_{j+1}$, $j\in \I_{n-1}$, $y\cdot e_n =  be_{1}$. 
\end{obs}

\subsection{Indecomposable modules}

Let $V \in \rep L$,  $\dim V = n \in \N$. 
Since $g^2$ is central, $V[\lambda]:=V_{g^2}^{(\lambda)} 
= V_{g}^{(a)}\oplus V_{g}^{(-a)} \leq V$  for any $\lambda = a^2 \in \ku^{\times}$. Let $\rep^{\lambda} L$ be the full subcategory
of $\rep L$ of those $V$ such that $V = V[\lambda]$. 

\begin{lema}\label{lema-generalized-g} 
\begin{enumerate}[leftmargin=*,label=\rm{(\roman*)}]
\item \label{it:repL-category}
$\rep L = \oplus_{\lambda\in \ku^{\times}} \rep^{\lambda} L$
is a graded tensor category, i.~e. $\otimes: \rep^{\lambda} L \times \rep^{\mu} L \to \rep^{\lambda\mu} L$ for  
$\lambda, \mu \in \ku^{\times}$, $\ku_{\varepsilon} \in \rep^{1} L$, the dual of $V \in  \rep^{\lambda} L$ belongs to
$\rep^{\lambda^{-1}} L$.
\item If $\Gamma \leq \ku^{\times}$, then $\rep^{\Gamma} L := \oplus_{\lambda\in \Gamma} \rep^{\lambda} L$
is a graded tensor subcategory of $\rep L$.
\item Let $\rep_1^{\lambda} L$ be the full subcategory
of $\rep^{\lambda} L$ consisting of modules with semisimple action of $g^2$.
Then $\rep^{\Gamma}_1 L := \oplus_{\lambda\in \Gamma} \rep_1^{\lambda} L$
is a graded tensor subcategory of $\rep L$ for every $\Gamma \leq \ku^{\times}$.
\item $\rep L^{(1)}$ can be identified with the tensor subcategory  $\rep_1^{1} L$. More generally,
$\rep L^{(n)}$ can be identified with $\rep_1^{\G_n} L$.
\item \label{it:comp-hom-repL} If $\lambda, b \in \ku^{\times}$, then $\leftaction_b, \rightaction_b: \rep^{\lambda} L  \to \rep^{\lambda b^2} L$,
$\leftaction_b(V) = \ku_{b} \otimes V$, $\rightaction_b(V) = V \otimes \ku_{b}$ are equivalences of abelian categories; 
in particular  $\leftaction_{-1}$ and $\rightaction_{-1}$ are auto-equivalences of  $\rep^{\lambda} L$ for all $\lambda$.
\item \label{it:indec-repL} If $V \in \rep L$ is indecomposable, then $V \in  \rep^{\lambda} L$ for a unique $\lambda$. 
\end{enumerate}
\end{lema} 
\pf
Let $\rep_n^{\lambda} L$ be the full subcategory
of $\rep L$ of those $V$ such that $V = \ker (g^2 - \lambda)^n$, $n \in \N_0$. We claim that 
$V \otimes W \in \rep_{nm}^{\lambda \mu} L$, whenever $V \in \rep_n^{\lambda} L$ and $W \in \rep_m^{\mu} L$.
The claim follows by induction on $\dim V$ and implies \ref{it:repL-category}. The rest of the proof is standard.
\epf

\begin{obs}
The various tensor categories $\rep^{\Gamma} L$  can be realized as the categories of comodules of a suitable Hopf algebra, namely
the Hopf subalgebra of the restricted dual of $L$ spanned by the matrix coefficients of the objects in $\rep^{\Gamma} L$.
Similarly for variations as in the Lemma.
\end{obs}

\begin{obs}
Since $\ku_a, \U_{a,b} \in \rep^{\lambda} L$ for any $\lambda = a^2 \in \ku^{\times}$, we see that
\begin{align*}
\ext_{L}^1(\ku_a, \ku_c) =  \ext_{L}^1(\ku_a, \U_{c,d}) = \ext_{L}^1(\U_{a,b}, \U_{c,d}) = 0,
\end{align*}
whenever $a^2 \neq c^2$.
\end{obs}

\subsubsection{Representations with $y = 0$}
Observe that $L/ LyL \simeq \ku \mathbb Z$. 
Thus there is a unique indecomposable module $\Vc_{a}^{n}$ of dimension $n$
where $y$ acts by $0$, namely with $g$ acting by $\jordan_n(a)$ with $a\in \ku^{\times}$.
 Given $n,m\in \N$ and $a\in \ku^{\times}$, it is easy to see that
$\Vc_{a}^n\simeq \ku_a\otimes \Vc_{1}^n$, $\Vc_{1}^n\otimes \Vc_{1}^m \simeq \Vc_{1}^m \otimes \Vc_{1}^n$.
Let  $n,m\in \N$, $2\leq m\leq n$. By \cite[Corollary 1]{Sri}, we have
$$\Vc_{1}^n\otimes \Vc_{1}^m \simeq \oplus_{k=1}^m\Vc_{1}^{n-(m-(2k-1))}.$$

\subsubsection{Representations with $y \neq 0$} We focus next on indecomposable $L$-modules with $y\neq 0$. 
Since $g$ has a unique eigenvalue $a$ if and only if $y$ acts by $0$,
 $g$ has eigenvalues $\pm a$. 
It would be enough to assume that $V \in \rep^{1} L$
since the indecomposable modules in $\rep^{\lambda} L$ can be deduced by Lemma \ref{lema-generalized-g} \ref{it:comp-hom-repL} but the analysis is the same as in the general case.
So assume that $V\in \rep^{\lambda} L$ indecomposable, $\lambda = a^2$.
Let $\base_{\pm}$ be a basis of $V_{g}^{(\pm a)}$ and $\base:=\base_+\cup \base_{-}$;
let $\ell = \dim V_{g}^{(a)}$, $\wp = n - \ell$. Up to replacing $V$ by $\leftaction_{-1}(V)$ 
(i.~e. interchanging $a$ and $-a$),
we may, and always will, assume that $\ell \geq \wp$. Then 
\begin{align} \label{matrixrep}
[g]_{\base} &=\begin{pmatrix}
A & 0 \\
0 &  B \end{pmatrix}, & 
[y]_{\base} &=\begin{pmatrix}
0 & C  \\
D & 0 \end{pmatrix},
\end{align}
where $A\in \GL_{\ell}(\ku)$, $B\in \GL_{\wp}(\ku)$, $C\in\operatorname{M}_{\ell\times \wp}(\ku)$ and $D\in\operatorname{M}_{\wp\times \ell}(\ku)$. 
Conversely,  $g$ and $y$ given by \eqref{matrixrep} define a representation of $L$ 
if and only if
\begin{align} \label{eqrep}
AC&=-CB, &  BD=-DA.&
\end{align}
We fix a pair of matrices $A$ and $B$ in Jordan form.
Our goals are to describe conditions on $y$ equivalent to the indecomposability of $V$ and then to classify the pairs
$(C,D)$ satisfying these conditions. In general, let 
\begin{align*}
\mathbb H &= \mathbb H_{A,B} = \{(t_1, t_2) \in \GL_{\ell}(\ku) \times \GL_{\wp}(\ku): t_1A = At_1, \, t_2B = Bt_2 \},
\\
\mathfrak V &= \mathfrak V_{A,B} = \{(C,D)\in\operatorname{M}_{\ell\times \wp}(\ku) \times\operatorname{M}_{\wp\times \ell}(\ku):
\text{\eqref{eqrep}  holds and $V$ is simple} \}.
\end{align*}
Then the algebraic group $\mathbb H$ acts on $\mathfrak V$ and our goals can be rephrased as:
\begin{itemize}[leftmargin=*]\renewcommand{\labelitemi}{$\circ$}
	\item Describe $\mathfrak V$ and decide when is non-empty,
	\item determine the orbits of the action of $\mathbb H$ on $\mathfrak V$.
\end{itemize}

We approach these questions by elementary means when $A$ and $B$ have simple Jordan forms
and apply the obtained results to the classification of the indecomposable modules of dimension $\leq 5$.

In this section, we set $U^{(\pm a)}_g=U^{(\pm a)}$ and $y_{\pm a} = y_{\vert U^{(\pm a)}}$
for any  $U \in \rep L$.

\subsubsection{Representations with $g$ semisimple}
Here $\dim V =n \geq 3$ and $g$ acts by $\begin{pmatrix}
a \id_{\ell} & 0 \\
0 &  -a \id_{\wp} \end{pmatrix}$ where $0 < \ell < n$ and $\wp=n-\ell$; any $C$, $D$ satisfy \eqref{eqrep}.

\begin{lema}\label{lema:g-semisimple-1} If either $y_{a} =0$ or $y_{-a} =0$, then $V$ is decomposable.
\end{lema}

\pf This is clear if $y_{a} =0$ and $y_{-a} =0$. If $y_{a} =0$ and $y_{-a} \neq 0$, then take a basis
$v_1, \dots, v_{q}$ of $\Ig(y_{-a})$, say $v_j = y_{-a}(v_{\ell + j})$, $j\in \I_q$, and a basis 
$v_{\ell + q + 1}, \dots, v_{\wp}$ of $\ker y_{-a}$.
Complete to a basis $v_1, \dots, v_{q}, \dots, v_{\ell}$ of $V^a$. 
Then $V = \langle v_1, v_{\ell + 1}\rangle \oplus  \langle v_2, \dots,  v_{\ell}, v_{\ell + 2}, \dots,  v_{\wp}\rangle$ is decomposable. If $y_{a}  \neq 0$ and $y_{-a} = 0$, then $\leftaction_{-1} (V)$ is decomposable by the preceding, and so is $V$.
\epf	

Assume that $\ell = n-1$. Define the representation
$\cC^n_{a}$ by $y\mapsto \begin{pmatrix}
0  & e_1 \\ 
e_{\ell}^t & 0  \end{pmatrix}$.

\begin{lema}\label{lema:g-semisimple-2} 
	Assume that $\ell = n-1$, $y_{a} \neq 0$  and $y_{-a} \neq 0$. 
	\begin{enumerate}[leftmargin=*,label=\rm{(\roman*)}]
		\item\label{it:2.5-1} If $y_{a}y_{-a} \neq 0$, then $V$ is decomposable.
		\item\label{it:2.5-2} $y_{a}y_{-a} = 0$ if and only if $V \simeq \cC^n_{a}$.
		\item\label{it:2.5-3} $\cC^n_{a}$ is indecomposable if and only if $n = 3$.
	\end{enumerate}
	
\end{lema}

\pf \ref{it:2.5-1} Assume that $y_{a}y_{-a} \neq 0$. Let $(v_{\ell + 1})$ be a basis of $V^{-a}$
and let $v_2, \dots, v_{\ell}$  be  a basis of $\ker y_{a}$. By hypothesis $v_1 := y_{-a}(v_{\ell + 1}) \notin \ker y_{a}$,
hence $v_1, \dots, v_{\ell}$ is a basis of $V^{a}$.  
Then $V = \langle v_1, v_{\ell + 1}\rangle \oplus \langle v_2, \dots,  v_{\ell}\rangle$ is decomposable.

\smallbreak
\noindent \ref{it:2.5-2} Since $n > 2$, $y_{a}y_{-a} = 0$ in $\cC^n_{a}$.
Assume that $y_{a}y_{-a} = 0$.  Let $(v_{\ell + 1})$ be a basis of $V^{-a}$ and
$v_1 := y_{-a}(v_{\ell + 1}) \in \ker y_{a}$; complete to a basis $v_1, \dots, v_{\ell - 1}$ of $\ker y_{a}$. 
Let  $v_\ell \in V^{a}$ be such that $y_{a}(v_\ell) = v_{\ell + 1}$. 
Then  $v_1, \dots, v_{\ell + 1}$ is a basis of $V$ that provides the isomorphism with  $\cC^n_{a}$.

\smallbreak
\noindent \ref{it:2.5-3} If $n > 3$, then $\langle e_1,  e_{\ell}, e_{\ell + 1}\rangle \oplus \langle e_2, \dots,  e_{\ell - 1}\rangle$ is a decomposition of $V = \cC^n_{a}$. Assume that $n = 3$ and that $V = U \oplus W$ is decomposable.
Hence $V^{\pm a} = U^{\pm a} \oplus W^{\pm a}$ so that $U^{-a} = 0$ and $V^{- a} =W^{-a} = \ku e_3$ or vice versa.
If $U^{-a} = 0$, then $U^{a} \leq \ker y_a = \ku e_1 = \ku  y_{-a}(e_3) \leq W^a$; hence $U =0$.
\epf

Assume that $\ell = n-2$. We define the representations $\cD_{1,a}^{4}$, $\cD_{2,a,b}^{4}$, $b\in \ku^{\times}$,
$\cD_{3,a}^{4}$ and $\cD_{4,a}^{5}$ by $y\mapsto \begin{pmatrix}
0 & C \\ D &  0 \end{pmatrix}$, where
\begin{align*}
&\cD_{1,a}^{4}:& & C=(0 \,\,\,e_1),& &D=\begin{pmatrix}
e_1^t \\
e_2^t  \end{pmatrix};&
&\cD_{2,a,b}^{4}:& & C=(e_1 \,\,\, e_2),& &D=\begin{pmatrix}
be_1^t+e_2^t \\
be_2^t  \end{pmatrix};
\\
&\cD_{3,a}^{4}:& & C=(e_1 \,\,\, e_2),& &D=\begin{pmatrix}
e_2^t \\
0  \end{pmatrix};&
&\cD_{4,a}^{5}:& & C=(e_1 \,\,\, e_2),& &D=\begin{pmatrix}
e_2^t \\
e_3^t  \end{pmatrix}.
\end{align*}

\begin{lema}\label{lema:g-semisimple-3}
 Let $V \in \rep L$ such that $n \geq 4$ and $g$ acts by $\begin{pmatrix}
	a \id_{\ell} & 0 \\
	0 &  -a \id_{2} \end{pmatrix}$. 
	\begin{enumerate}[leftmargin=*,label=\rm{(\roman*)}]
\item\label{item:g-semisimple-3-n4} If $n=4$, then $V$ is indecomposable iff $V \simeq  \cD_{1,a}^{4}$,
or $\cD_{2,a,b}^{4}$ or $\cD_{3,a}^{4}$.
\item\label{item:g-semisimple-3-n5} If $n=5$, then $V$ is indecomposable iff  $V \simeq \cD_{4,a}^{5}$.
\item\label{item:g-semisimple-3-n6} If $n \geq 6$, then $V$ is decomposable.
	\end{enumerate}	
\end{lema}

\pf
\begin{step}
The $L$-modules $\cD_{1,a}^{4}$, $\cD_{2,a,b}^{4}$, $\cD_{3,a}^{4}$ and $\cD_{4,a}^{5}$ are indecomposable.
\end{step}
 
Suppose that $V = \cD_{1,a}^{4}$ is decomposable, i. e. there exist non-zero $L$-submodules
$U$ and $W$  such that $V=W\oplus U$. Then $V^{-a}=W^{-a}\oplus U^{-a}$. If $W^{-a} =0$, then $yW^{a} =0$,
hence $W^{a} =0$ since $y_{a}$ is injective. Thus $\dim W^{-a} = \dim U^{-a} = 1$.
If $0\neq w=\alpha e_3+\beta e_4 \in W^{-a}$, then $y^2(w)= \beta e_3 \in  W^{-a}$, so that either $\alpha =0$ or 
$\beta =0$, and similarly for $U^{-a}$. In other words, we may assume that 
$W^{-a} = \la e_3 \ra$ and $U^{-a} = \la e_4 \ra$, but $y^2(e_4)=  e_3 \in  W^{-a} \cap U^{-a}$, a contradiction.
By a similar argument, $\cD_{4,a}^{5}$ is indecomposable.
Clearly $\cD_{3,a}^{4} \simeq \left(\cD_{1,a^{-1}}^{4} \right)^{*}$, consequently it is indecomposable.

Finally, suppose that $\cD_{2,a,b}^{4}=U\oplus W$ for some $L$-submodules $U$ and $W$.
Pick a square root $\sqrt{b}$ of $b$ and note that $\lambda_i=(-1)^i\sqrt{b}$, $i=1,2$, are the eigenvalues of $y$ and $V_{y}^{\lambda_i}=\langle v_i:=e_1+ \lambda_i e_3 \rangle$ are the respective eigenspaces.  
If $v_1=u+w$ with $u\in U$ and $w\in W$, then $u, w \in V_{y}^{\lambda_1}$. If $u\neq 0$ and $w\neq 0$, then $v_1\in U\cap W= 0$, a contradiction. If $w=0$, then $v_1\in U$ and $v_2=gv_1=gu\in U$ which implies $W=0$. Similarly, $u=0$ implies $U=0$.

\begin{step}\label{aux1}
	If $n \geq 5$ and $y_ay_{-a}$ is an isomorphism, then $V$ is decomposable.
\end{step}

Note that $y_ay_{-a}$ is an isomorphism if and only if $y_{-a}$ is injective, $y_a$ is surjective and $\Ig y_{-a} \cap \nuc y_a = 0$. Let $e_{\ell +1}, e_{\ell+2}$ be a basis of $V^{-a}=\Ig y_{a}$. Denote by $e_i=y_{-a}(e_{i+2}) \neq 0$, where $i \in \I_{\ell-1,\ell}$. So $e_i \notin \nuc y_a$. We complete to a basis $e_1, \cdots, e_{\ell}$ of $V^a$ such that $y_a(e_j)=0$, for all $j \in \I_{1,\ell-2}$. Therefore $V= \la e_1, \cdots, e_{\ell-2} \ra \oplus \la e_{\ell-1}, e_{\ell}, e_{\ell+1}, e_{\ell +2}\ra$ is decomposable.

\smallbreak
We next investigate what happens when $y_ay_{-a}$ is not an isomorphism. 

\begin{step}\label{aux2}
	Assume that $y_{-a}$ is injective and $\Ig y_{-a} \cap \nuc y_a \neq 0$.
	\begin{enumerate}[leftmargin=*,label=\rm{(\alph*)}]
\item\label{item:aux2-a} If $n \geq 6$ and $y_{a}$ is surjective then $V$ is decomposable.
\item\label{item:aux2-b}  If $n \geq 5$ and $y_{a}$ is not surjective then $V$ is decomposable.
	\end{enumerate}	
\end{step}

Let $0 \neq e_1 \in \Ig y_{-a} \cap \nuc y_a$. 
Pick $e_{\ell+1} \in V^{-a}$ such that $e_1=y_{-a}(e_{\ell+1})$ and
complete to a basis $e_{\ell+1}, e_{\ell+2}$ of $V^{-a}$. 
Let $e_2 := y_{-a}(e_{\ell+2})$.

\ref{item:aux2-a}:  If $e_2 \in \nuc y_a$, then 
complete to a basis $e_1, \dots, e_{\ell-2}$ of $\nuc y_a$. 
Take $e_{\ell-1}, e_{\ell} \in V^{a}$ such that $y_a(e_{\ell-1})=e_{\ell+1}$ 
and $y_a(e_{\ell})=e_{\ell+2}$. Clearly $e_1, \dots, e_{\ell}$ is a basis of $V^a$. 
Since $\ell \geq 4$, 
$V= \la e_1,e_{\ell-1}, e_{\ell+1} \ra \oplus \la e_2, \cdots, e_{\ell-2}, e_{\ell}, e_{\ell+2} \ra$
is decomposable. 
If $e_2 \notin \nuc y_a$, then consider $e_1,e_3,\cdots , e_{\ell-1}$ a basis of  $\nuc y_a$. 
We can find $e_{\ell} \in V^{a}$ 
such that $e_1,e_3, \cdots, e_{\ell-1}, e_2, e_{\ell}$ is a basis of $V^{a}$. 
Since $\ell \geq 4$, $V=\la e_1,e_2,e_{\ell}, e_{\ell+1},e_{\ell+2}\ra \oplus \la e_3, \cdots,e_{\ell-1} \ra$ 
is decomposable.

\ref{item:aux2-b}: If $e_2 \in \nuc y_a$, then we complete to a basis $e_1,e_2 \cdots, e_{\ell-1}$ 
of $\nuc y_a$ and then to a basis $e_1, \cdots, e_{\ell}$ of $V^a$. 
So $y_a(e_{\ell})=be_{\ell+1}+ce_{\ell+2}$, $b,c \in \ku,\, bc \neq 0$. 
Then either $V= \la be_1+ce_2,e_{\ell},be_{\ell+1}+ce_{\ell+2} \ra \oplus \la e_2,\cdots, e_{\ell-1}, e_{\ell+2}\ra$
if $b \in \ku^{\times}$, or else $V= \la e_1, e_3, \cdots,e_{\ell-1}, e_{\ell+1} \ra \oplus \la e_2,e_{\ell}, e_{\ell+2} \ra$
if $b=0$ and $c \in \ku^{\times}$. That is, $V$ is decomposable. 
If $e_2 \notin \nuc y_a$, take a basis $e_1,e_3,\cdots , e_{\ell}$ of  $\nuc y_a$, so that $e_1,e_2, \cdots, e_{\ell}$ is a basis of $V^{a}$. Then $V=\la e_1,e_2, e_{\ell+1},e_{\ell+2}\ra \oplus \la e_3, \cdots,e_{\ell} \ra$ is decomposable.

\begin{step}\label{aux3}
If $n \geq 4$, $y_{-a}$ is not injective and $\Ig y_{-a} \cap \nuc y_a \neq 0$, then $V$ is decomposable.
\end{step}

 Let $0 \neq e_1 \in \Ig y_{-a} \cap \nuc y_a$. Pick $e_{\ell+2} \in V^{-a}$ such that $e_1=y_{-a}(e_{\ell+2})$ and complete to a basis $e_{\ell+1}, e_{\ell+2}$ of $V^{-a}$ with $y_{-a}(e_{\ell+1})=0$. If $y_a$ is surjective, take $e_{\ell+1}=y_a({e_{\ell-1}})$ and $e_{\ell+2}=y_a({e_{\ell}})$ and complete to a basis $e_1,e_2,\dots, e_{\ell}$ of $V^a$ such that $e_i \in \nuc y_a$, $i \in \I_{\ell - 2}$. Then $V= \la e_1, e_{\ell},e_{\ell+2} \ra \oplus \la e_2, \dots, e_{\ell-1}, e_{\ell+1}\ra$ is decomposable. If $y_{a}$ is not surjective, then take $e_1,e_2 \dots, e_{\ell-1}$ be a basis of $\nuc y_a$ such that  $e_1,e_2 \dots, e_{\ell}$ is a basis of $V^a$. Then $V=  \la e_1,e_3, \dots, e_{\ell+2} \ra \oplus \la e_2\ra$ is decomposable.

\begin{step}\label{aux4}
	Assume that  $y_{-a}$ is not injective and $\Ig y_{-a} \cap \nuc y_a = 0$.
\begin{enumerate}[leftmargin=*,label=\rm{(\alph*)}]
\item\label{item:aux4-a} If $n \geq 5$ and $y_{a}$ is surjective then $V$ is decomposable.
\item\label{item:aux4-b}  If $n \geq 4$ and $y_{a}$ is not surjective then $V$ is decomposable.
\end{enumerate}	
\end{step}

\ref{item:aux4-a}:  Let $e_{\ell+1},e_{\ell+2}$ be a basis of $V^{-a}$ such that $y_{-a}(e_{\ell+1})=0$ and $y_{-a}(e_{\ell+2})$ $= e_{\ell} \neq 0$, consequently $e_{\ell} \notin \nuc y_a$. Thus we can consider $e_1,e_2, \dots, e_{\ell}$ a basis of $V^a$ such that $e_1,e_2, \dots, e_{\ell-2} \in \nuc y_a$. Therefore $V= \la e_1, \dots, e_{\ell-2}\ra \oplus \la e_{\ell-1}, e_{\ell},e_{\ell+1},e_{\ell+2} \ra$ is decomposable.
\ref{item:aux4-b} is similar to \ref{item:aux4-a}.

\begin{step}\label{aux5}
Proof of the Lemma.
\end{step}

\ref{item:g-semisimple-3-n4} By Steps \ref{aux1}, \ref{aux2}, \ref{aux3} and \ref{aux4} it is enough to consider the cases:
\begin{itemize}[leftmargin=*]\renewcommand{\labelitemi}{$\circ$}
\item  $y_ay_{-a}$ is an isomorphism.
\item $ y_{-a}$ is injective and $ y_{a}$ is not surjective.
\item $y_{-a}$ is not injective and $ y_{a}$ is surjective. 
\end{itemize} 
 
In the first case we have two possibilities for the canonical Jordan form of $y_ay_{-a}$, namely diagonal or $J_2(b)$, with $b \in \ku^{\times}$. Consider $e_3,e_4$ basis of $V^{-a}$ such that the canonical Jordan form of $y_ay_{-a}$ is diagonal. Take $e_1=y(e_3)$ and $e_2=y(e_4)$. Then $V= \la e_1,e_3 \ra \oplus \la e_2, e_4 \ra$ is decomposable. Now if the canonical Jordan form of $y_ay_{-a}$ is $J_2(b)$ then $V= \cD_{2,a,b}^{4}$. 

In the second case, if $y_ay_{-a}$ is diagonalizable then $V$ is decomposable and otherwise $V=\cD_{3,a}^{4}$. 

Finally in the third case, consider $e_3 \in \nuc{y_{-a}}$ and complete to a basis $e_3,e_4$ of $V^{-a}$ such that $y_ay_{-a}(e_4)=be_4$, $b \in \ku^{\times}$. Take $e_1=y_{-a}(e_4)$ and complete to a basis $e_1,e_2$ of $V^a$ with $y_a(e_2)=ce_3+de_4$, $c \in \ku^{\times}$. So $V=\la e_1,e_4\ra \oplus \la -db^{-1}e_1+e_2, e_3\ra $ is decomposable. Consider $e_3,e_4$ a basis of $V^{-a}$ such that the canonical Jordan form of $y_ay_{-a}$ is $J_2(0)$. Take $e_1=y_{a}(e_4)$, $e_1,e_2$ a basis of $V^a$ with $y_a(e_2)=be_3+ce_4$, $c \in \ku^{\times}$. Taking the basis $e_1, c^{-1}e_2,e_3, bc^{-1}e_3+e_{4}$ we have that $V \simeq \cD_{1,a}^{4}$.

\ref{item:g-semisimple-3-n5} By Steps \ref{aux1}, \ref{aux2}, \ref{aux3} and \ref{aux4} it is enough to consider the case where $y_{-a}$ is injective, $y_{a}$ is surjective and $\Ig y_{-a} \cap \nuc y_a \neq 0$. Let $0 \neq e_1 \in \Ig y_{-a} \cap \nuc y_a$. Pick $e_{4} \in V^{-a}$ such that $e_1=y_{-a}(e_{4})$  and complete to a basis $e_{4}, e_{5}$ 
of $V^{-a}$. Denote by $e_2=y_{-a}(e_{5})$. Let $ye_2=\alpha e_4+\beta e_5$, $ye_3= \gamma e_4 + \eta e_5$, $A=\left(
\begin{array}{cc}
\alpha & \gamma \\
\beta & \eta \\
\end{array}
\right)$. Taking the basis $e_1, \eta (\dete A)^{-1}e_2-\beta (\dete A)^{-1}e_3, -\gamma (\dete A)^{-1}e_2+\alpha (\dete A)^{-1}e_3, e_4,e_5$, we have that $V \simeq \cD_{4,a}^{5}$.

\ref{item:g-semisimple-3-n6}  follows from the Steps \ref{aux1}, \ref{aux2}, \ref{aux3} and \ref{aux4}.
\epf

\subsubsection{Representations with $A$ a Jordan block}

\begin{lema}\label{lema:two-blocks} Let $V$ be a representation of $L$ given by \eqref{matrixrep} where
$A$ and $B$ are Jordan blocks. Then $V$ is indecomposable if and only if $y \neq 0$.
\end{lema}

\pf Let $V = U \oplus W$ be a decomposition with both $U$ and $W$ not 0.
Hence $V^{(\pm a)} = U^{(\pm a)} \oplus W^{(\pm a)}$.
Since $A$ is a Jordan block, either $U^{(a)} = 0$ and $W^{(a)} =V^{(a)}$ or vice versa.
If $U^{(a)} = 0$, then $U^{(-a)} = V^{(-a)}$ and $W^{(-a)} =0$, because $B$ is also a Jordan block; 
thus $y_{-a} = yU^{(-a)} \leq U^{(a)} = 0$ and similarly $y(V^{(a)}) = 0$, implying
$y=0$. Same if $U^{(a)} = V^{(a)}$. The converse is evident.
\epf

Assume that $n \geq 3$ and $\ell = n-1$. Define the representations
$\Ec^n_{1,a}$, $\Ec^n_{2,a}$ and $\Ec^n_{3,a,b}$, $b \in \ku^{\times}$ by 
\eqref{matrixrep}, where $A = \jordan_{\ell}(a)$, $B = -a$,
and $y$ acts as follows: 
\begin{align*}
&\Ec^n_{1,a}: y\mapsto \begin{pmatrix} 0  & e_1  \\ 0  & 0 \end{pmatrix},
&\Ec^n_{2,a}: y\mapsto \begin{pmatrix} 0  & 0  \\ e_{\ell}^t  & 0 \end{pmatrix},&
&\Ec^n_{3,a,b}: y\mapsto \begin{pmatrix} 0  & be_1  \\  e_{\ell}^t  & 0 \end{pmatrix}.
\end{align*}
These  are pairwise non-isomorphic, for different values of $a$ and $b$.

\begin{lema}\label{lema:two-blocks-p=1} Let $n \geq 3$ and $\ell = n-1$.
Let $V$ be a representation of $L$ given by \eqref{matrixrep} where
$A =  \jordan_{\ell}(a)$ and $B = -a$. 
Then the following are equivalent:
\begin{enumerate}[leftmargin=*,label=\rm{(\roman*)}]
\item\label{item:two-blocks-p=1-1} $V$ is indecomposable,
\item\label{item:two-blocks-p=1-2} $y \neq 0$,
\item\label{item:two-blocks-p=1-3} $V$ is isomorphic to one, and only one, of $\Ec^n_{1,a}$,$\Ec^n_{2,a}$ or $\Ec^n_{3,a,b}$, for unique $a,b \in \ku^{\times}$.
\end{enumerate}
\end{lema}
 
\pf Here $C = v$, $D = w^t$ where $v, w\in \ku^{\ell}$ and \eqref{eqrep} says that $Av = av$, $A^tw = aw$. 
Thus $v = be_1$, $w = ce_{\ell}$ for some $b,c \in \ku$, 
and $V$ is indecomposable iff $(b, c) \neq 0$ by Lemma \ref{lema:two-blocks}.
If $c=0$, then the basis $be_1, \dots, e_n$  gives $V \simeq \Ec^n_{1,a}$. 
If $c \neq 0$, then  $e_1, \dots, ce_n$ gives $V \simeq \Ec^n_{2,a}$ when $b=0$, or 
$\Ec^n_{3,a,bc}$  when $b \neq 0$. 
\epf 

Assume that $n \geq 4$ and $\ell = n-2$, so that $A = \jordan_{\ell}(a)$ and $B =\jordan_{2}(-a)$. 
Consider $C_1,C_2\in\operatorname{M}_{\ell\times \wp}(\ku) $ and $D,D_i\in\operatorname{M}_{\wp\times\ell}(\ku)$, given by:
\begin{align*} 
C_1 &=  (0  \,\,\,\,\, e_1),& 
C_2 &=  (e_1  \,\,\,\,\, -e_2),& 
D &=  \begin{pmatrix} e^t_l  \\ 0 \end{pmatrix}, &
D_i &=  \begin{pmatrix}  	e_{i }^t  \\ -e^t_{i+1} \end{pmatrix},& i\in \I_{1,\ell-1}.
\end{align*}
We define the following representations of $L$ on $\ku^{n}$ by \eqref{matrixrep} where $y$ acts by:

\begin{align*}
&\Fc^n_{1,\,a}: y\mapsto \begin{pmatrix} 0  & C_1  \\ 0  & 0 \end{pmatrix},
&&\Fc^n_{2,\,a}: y\mapsto \begin{pmatrix} 0  & C_2  \\ 0  & 0 \end{pmatrix},&\\
&\Fc^n_{3,\,a}: y\mapsto \begin{pmatrix} 0  & 0  \\ D  & 0 \end{pmatrix},
&&\Fc^{n,i}_{4,\,a}: y\mapsto \begin{pmatrix} 0  & 0  \\ D_i  & 0 \end{pmatrix},&\\
&\Fc^n_{5,\,a,\,b,\,c}: y\mapsto \begin{pmatrix} 0  & bC_1  \\ cD  & 0 \end{pmatrix},&
&\Fc^{n,i}_{6,\,a,\,b,\, c}: y\mapsto \begin{pmatrix} 0  & bC_1  \\ cD_i  & 0 \end{pmatrix},&\\
&\Fc^n_{7,\,a,\,b,\,c}: y\mapsto \begin{pmatrix} 0  & bC_2  \\ cD  & 0 \end{pmatrix},&
&\Fc^{n,i}_{8,\,a,\,b,\,c}: y\mapsto \begin{pmatrix} 0  & bC_2  \\ cD_i  & 0 \end{pmatrix}.
\end{align*}

\begin{lema}\label{lema:two-blocks-p=2} 
	Let $V$ be a representation of $L$ given by \eqref{matrixrep} where
	$A =  \jordan_{\ell}(a)$ and $B = \jordan_{2}(-a) $, with $n \ge 4$. 
	Then the following are equivalent:
	\begin{enumerate}[leftmargin=*,label=\rm{(\roman*)}]
		\item\label{item:two-blocks-p=2-1} $V$ is indecomposable.

		\item\label{item:two-blocks-p=2-3} $V$ is isomorphic to one of $\Fc^n_{j,\,a}$, $j \in \I_{3}$, $\Fc^{n,i}_{4,\,a}$, $\Fc^n_{j,\,a,\,b,\,c}$, $j=5,7$ and $\Fc^{n,i}_{j,\,a,\,b,\,c}$, $j=6,8$, for unique $a,b,c \in \ku^{\times}$.
	\end{enumerate}
If $n \geq 5$, then the above representations are pairwise non-isomorphic.
If $n=4$, then $\Fc^{4,1}_{4,\,a} \simeq \Fc^{4}_{3,\,-a}$, $\Fc^{4}_{2,\,a} \simeq \Fc^{4,1}_{4,\,-a}$, $\Fc^{4}_{7,\,a,\,b,\,c} \simeq \Fc^{4,1}_{6,\,-a,\,c\,b}$.

\end{lema}
\pf By Lemma \ref{lema:two-blocks}, we may assume that $y \neq 0$. Then the Lemma follows by a lengthy but straightforward 
analysis. \epf 

Let $n \geq 4$, $\ell = n-2$ and $\wp = 2$.
Define the representation $\cG^n_{a}$ by 
\eqref{matrixrep}, where $A  = \jordan_{\ell}(a)$, $B = -a\Id_2$,
$C = (0 \,  e_1)$ and $D = {\small \begin{pmatrix}  e_{\ell}^t \\ 0	\end{pmatrix}}$.

\begin{lema}\label{lema:A-block-B-id-p>1} Assume that $\ell, \wp \ge 2$. 
Let $V$ be a representation of $L$ given by \eqref{matrixrep} where
$A =  \jordan_{\ell}(a)$ and $B = -a \Id_{\wp}$.  
Then the following are equivalent:
\begin{enumerate}[leftmargin=*,label=\rm{(\roman*)}]
	\item $V$ is indecomposable.
	
	\item $\wp = 2$ and $V \simeq \cG^n_{a}$.
\end{enumerate}
\end{lema}

\pf  Let $C\in\operatorname{M}_{\ell\times \wp}(\ku)$ and $D\in\operatorname{M}_{\wp\times \ell}(\ku)$ such that \eqref{eqrep} holds.
By our hypotheses on $A$ and $B$, 
there exist $c_1, \dots, c_{\wp}$, $d_1, \dots, d_{\wp} \in \ku$
such that $$ C=\sum_{i\in \I_{\wp}}c_ie_i,\qquad\quad  D=\sum_{i\in \I_{\wp}}d_ie^t_i.$$
Let $y_{\pm a} = y_{\vert V^{(\pm a)}}$. Since $\dim \Ig y _{(\pm a)} \leq 1$ and $\ell, \wp \geq 2$, we see that $\ker y_{(\pm a)} \neq 0$.

Assume first that $y_{a} = 0$. Let $U = U^{(a)} \oplus U^{(-a)}$, where $U^{(a)} = 0$ and $U^{(-a)} = \ker y_{-a} \neq 0$;
let  $W = W^{(a)} \oplus W^{(-a)}$, where $W^{(-a)}$ is a direct summand of $\ker y_{-a}$ in $V^{(-a)}$ and $W^{(a)} = V^{(a)}$. Then 
$V = U \oplus W$ is a decomposition in $\rep L$.
So, assume that $y_{a} \neq 0$. 

\begin{itemize}[leftmargin=*]\renewcommand{\labelitemi}{$\circ$}
	
\item If $\Ig y_{a} \cap \ker y_{-a} =0$, then pick 
a direct summand $Z$ of $\Ig y_{a} \oplus \ker y_{-a}$ in $V^{(-a)}$ and set $U =  U^{(-a)} = \ker y_{-a}$, 
$W = V^{(a)} \oplus (\Ig y_{a} \oplus Z)$.

\item If $\Ig y_{a} \lneq \ker y_{-a}$, then pick 
a direct summand $Z_1$ of $\Ig y_{a}$ in $\ker y_{-a}$ and a direct summand $Z_2$ of $\ker y_{-a}$ in $V^{(-a)}$.
Then set $U =  U^{(-a)} = Z_1 \neq 0$, 
$W = V^{(a)} \oplus (\Ig y_{a} \oplus Z_2)$.
\end{itemize}
Hence $V = U \oplus W$ is a decomposition in $\rep L$ in both cases.

It remains the case $\Ig y_{a} = \ker y_{-a}$, necessarily $\wp = 2$ and $y_{-a} \neq 0$.
Let $v_1, v_2 \in \ku^2$ be such that $\Ig y_{a} = \ku v_1$ and  $y_{-a}(v_2) =  e_1$.
Then $\ku^2 = \ku v_1 \oplus \ku v_2$ and considering the basis $e_1, \dots, e:{\ell}, v_1, v_2$, we see that
$V \simeq \cG^n_{a}$.

Finally, we show that $V = \cG^n_{a}$ is indecomposable. Indeed, let 
$V = U \oplus W$ be a decomposition in $\rep L$. Then $V^{(a)} = U^{(a)} \oplus W^{(a)}$ and either 
$U^{(a)}=0$ or $W^{(a)} = 0$; suppose the first happens. Then $U^{(-a)} \leq \ker y_{-a} = y_{a}(V^{(a)}) = 
y_{a}(W^{(a)}) \leq W^{(-a)}$, hence $U^{(-a)} =0$ and a fortiori $U =0$. 
\epf

\subsubsection{Representations with $A$ sum of a Jordan block and a line} 
Assume $n\geq 4$.
We fix $A = \begin{pmatrix}\jordan_{n-2}(a) & \\ & a \end{pmatrix}$, $B = -a$.
Define the representations $\cH_{1,a}^{n}$, $\cH_{2,a}^{n}$, $\cH_{3,a,b}^{n}$, $b\in \ku^{\times}$ by \eqref{matrixrep},  
where $y$ acts on $\ku^n$ by:
\begin{align*} \scriptscriptstyle
&\cH_{1,a}^{n}: y \mapsto{\small \begin{pmatrix}
0  & e_1 \\
e_{n-1}^t & 0  \end{pmatrix}}, \,
\cH_{2,a}^{n}:  y \mapsto {\small \begin{pmatrix}
0  & e_{n-1} \\
e_{n-2}^t& 0 \end{pmatrix}}, \,
\cH_{3,a,b}^{n}:  y \mapsto {\small \begin{pmatrix}
0  & e_{n-1} \\
be_{n-1}^t& 0 \end{pmatrix}}.
\end{align*}

It is easy to see that these modules are indecomposable.

\begin{lema}\label{lema:jordan-plus-1dim}  Let $V$ be a representation of $L$ on $\ku^n$  given by \eqref{matrixrep}, with $A,B$ as above. Then the following statements are equivalent:
	\begin{enumerate}[leftmargin=*,label=\rm{(\roman*)}]
		\item $V$ is indecomposable.
		\item $V$ is isomorphic to  exactly one of $\cH_{1,a}^{n}$, $\cH_{2,a}^{n}$ or $\cH_{3,a,b}^{n}$, for unique $b\in \ku^{\times}$.
	\end{enumerate}		
\end{lema}

\pf 
\begin{passo} There are $b,c,d,f \in \ku$ such that
$C=be_1+ce_{n-1}$, $D=de_{n-2}^{t}+fe_{n-1}^{t}$.
\end{passo}

This follows at once from \eqref{eqrep}.

\begin{passo}\label{passo:dos} If either $y_{-a}=0$ or $y_{a}=0$ then $V$ is decomposable.
\end{passo}

Assume $y_{-a}=0$, i.~e. $b=c=0$. If $d=0$, then $\left\langle e_1,\ldots,e_{n-2} \right\rangle \oplus \left\langle e_{n-1}, e_n \right\rangle$ is a decomposition of $V$. If $f=0$, then $\left\langle e_1,\ldots,e_{n-2},e_n \right\rangle \oplus \left\langle e_{n-1} \right\rangle$ is a decomposition of $V$. Also  $\left\langle d^{-1}e_1,\ldots,d^{-1}e_{n-3},d^{-1}e_{n-2}-f^{-1}e_{n-1}\right\rangle \oplus \left\langle e_{n-1},e_n \right\rangle$ is a decomposition of $V$ when $d,f\neq 0$. The case $y_{a}=0$ is similar.

\begin{passo}\label{passo:tres}  $V$ is indecomposable iff
	$(be,cd)\neq 0$.
\end{passo}

Suppose that $(be,cd)=0$. If $b=0$, then $c\neq 0$ by Step \ref{passo:dos}. 
Thus $d=0$ and $\left\langle e_1,\ldots,e_{n-2} \right\rangle \oplus \left\langle e_{n-1}, e_n \right\rangle$ is a decomposition of $V$. Similarly, if $e=0$, then  $d\neq 0$,  $c=0$ and $\left\langle e_1,\ldots,e_{n-2},e_n \right\rangle \oplus \left\langle e_{n-1} \right\rangle$ is a decomposition of $V$.

\smallbreak
Conversely, suppose that  $V =U\oplus W$ with $U,W$ non-trivial submodules of $V$. 
From $V^{(-a)}=U^{(-a)}\oplus W^{(-a)}$, we may assume that $W^{(-a)}=0$ and $U^{(-a)}=V^{(-a)}$. Since $V^{(a)}=U^{(a)}\oplus W^{(a)}$
and $yW^{(a)}\subset W^{(-a)}=0$,
there are three possibilities:
if $W^{(a)} = \langle e_{n-1}\rangle$, then $f=0$;
if $W^{(a)} = \la e_1,\ldots,e_{n-2} \ra$, then $d=0$; if 
$W^{(a)} = V^{(a)}$, then $d = f =0$.
	
\begin{passo} If $V$ is indecomposable,  then it is one of $\cH_{1,a}^{n}$, $\cH_{2,a}^{n}$ or $\cH_{3,a,b}^{n}$.
\end{passo}

By Step \ref{passo:tres}, $(bf,cd)\neq 0$ and we proceed by a straightforward analysis. \epf

\subsubsection{Representations with $A$ a sum of Jordan blocks of the same size}
Let  $r, t \ge 2$ such that $\ell = rt$ and $\wp \geq 1$. Assume that $A$ consists of $t$ blocks of size $r$, i~e.
$A = {\small \begin{pmatrix}
	\jordan_{r}(a) & & \\ & \ddots & \\ & &\jordan_{r}(a)
	\end{pmatrix}}$. Let $y_{\pm a} = y_{\vert V^{(\pm a})}$. 

Assume that $t = 2$, thus $\ell = 2r$ and $\wp =1$. Define the representation
$\mathcal{I}^n_{a}$ by $y\mapsto \begin{pmatrix}
0  & e_1 \\ 
e_{2r}^t & 0  \end{pmatrix}$.

\begin{lema}\label{lema:A-t-bloques-r-B-id-p1} Assume that $\wp =1$. 
	Let $V$ be a representation of $L$ given by \eqref{matrixrep} where
	$A$ is as above and $B = -a$.  
	Then the following are equivalent:
	\begin{enumerate}[leftmargin=*,label=\rm{(\roman*)}]
		\item $V$ is indecomposable,
		
		\item $y_{a} \neq 0$,  $y_{-a} \neq 0$, $y_{a}y_{-a} = 0$ and $t = 2$,
		
		\item  $V\simeq \mathcal{I}^n_{a}$.
	\end{enumerate}
\end{lema}

\pf Let $C\in\operatorname{M}_{\ell\times \wp}(\ku)$ and $D\in\operatorname{M}_{\wp\times \ell}(\ku)$ such that \eqref{eqrep} holds.
By our hypotheses on $A$ and $B$, 
there exist $c_1, \dots, c_{t}$, $d_1, \dots, d_{t} \in \ku$
such that 
\begin{align*}
C &= \sum_{j\in \I_t} c_j e_{(j-1)r + 1} ,&
D &= \sum_{j\in \I_t} d_j e^t_{jr} .
\end{align*}
We summarize some well-known facts about the $g$-submodules of $V^{a}$: 
\begin{enumerate}[leftmargin=*,label=\rm{(\alph*)}]
\item Let $T = \langle e_{jr}: j\in \I_t\rangle$. If 
$w = \sum_{j\in \I_t} \alpha_j e_{jr} \in T$, then
\begin{align}\label{eq:base-wg}
\langle w\rangle_g = \big\langle \sum_{j\in \I_{t}} \alpha_{j} e_{jr +1 -i}: i \in  \I_{r}\big\rangle.
\end{align} 
\item  If  $T = R \oplus S$, then 
$V^{(a)} = \langle R\rangle_g\oplus \langle S\rangle_g$.
\end{enumerate}

Clearly, $V$ is decomposable if $y_{a} =0$ and $y_{-a} =0$.
Observe that if $y_{-a} \neq 0$, then $\Ig y_{-a} \langle C\rangle \subseteq \langle w_1\rangle_g$ 
	where $w_1 = \sum_{j\in \I_{t}} c_{j} e_{jr}$.

Assume that $y_{a} =0$ and $y_{-a} \neq 0$. Then complete $w_1$ to a basis  $w_1, \dots, w_t$ of $T$. 
Set $U = \oplus_{j\in \I_{2,t}} \langle w_j\rangle_g$, $W = \langle w_1\rangle_g \oplus V^{-a}$.
Then $V = U \oplus W$ is a decomposition in $\rep L$.

Assume that $y_{a} \neq 0$ and $y_{-a} = 0$. If $w \in T \cap \ker y_{a}$, 
then $\langle w\rangle_g  \subseteq \ker y_{a}$. By our present hypothesis, $T \cap \ker y_{a} \neq T$.
Pick $w \in T - T \cap \ker y_{a}$ and set $U = \langle T \cap \ker y_{a}\rangle_g$, 
$W = \langle w\rangle_g \oplus V^{(-a)}$.
Then $V = U \oplus W$ is a decomposition.

Assume that $y_{a} \neq 0$ and $y_{-a} \neq 0$. If $y_{a}y_{-a} \neq 0$, then $w_1 \notin T \cap \ker y_{a}$.
Then $V = \langle T \cap \ker y_{a}\rangle_g \oplus (\langle w_1\rangle_g \oplus V^{(-a)})$ is a decomposition in $\rep L$.
If otherwise $y_{a}y_{-a} = 0$, then $w_1 \in T \cap \ker y_{a}$; 
complete  to a basis  $w_1, \dots, w_{t-1}$ of $T \cap \ker y_{a}$ and pick $w_t \in T$: $y_{a}(w_t) = 1$.
Take in this case $U = \langle w_2, \dots, w_{t-1} \rangle_g$, 
$W = \langle w_1, w_t\rangle_g \oplus V^{(-a)}$.
Then $V = U \oplus W$ is a decomposition in $\rep L$ unless $t=2$ in which case $U = 0$.
When $t = 2$ then the basis of $V^{(a)}$ obtained from those of $\langle w_1\rangle_g$ and  
$\langle w_2 \rangle_g$ given by \eqref{eq:base-wg} realizes the isomorphism with 
$\mathcal{I}^n_{a}$.

Finally we show that $\mathcal{I}^n_{a}$ is indecomposable. 
Let $V = U \oplus W$ be a decomposition. Then $V^{(\pm a)} = U^{(\pm a)} \oplus W^{(\pm a)}$.
We may assume that $V^{(-a)} = U^{(-a)}$, $W^{(-a)} = 0$. Then
$W^{(a)} \subseteq \ker y_{a} = \Ig y_{-a} \subseteq U^{(a)}$, hence $W =0$.
\epf

\subsection{Low dimension}
In this subsection we classify those $V \in \rep L$ indecomposable of dimension $n\leq 5$.
If $g$ has a unique eigenvalue $a$, then $V \simeq \Vc_{a}^n$. 
So in all proofs below, $g$ has two eigenvalues $\pm a$ and the representation is given by matrices as in 
\eqref{matrixrep} satisfying \eqref{eqrep}. Also we assume that  
$\dim V^{(a)} =\ell \geq \wp = \dim V^{(-a)}$, otherwise apply  $\leftaction_{-1}$.

\subsubsection{Dimension 2}
Given $a\in \ku^{\times}$, let $\W_{a} = \ku^2$ be the representation given by
$g \mapsto \begin{pmatrix}
a & 0 \\
0 & -a \end{pmatrix}$,  $y\mapsto \begin{pmatrix}
0 & 1 \\
0 & 0 \end{pmatrix}$.
By Lemma \ref{lema:two-blocks}, $\W_{a}$ is indecomposable.

\begin{prop}\label{prop:L-dim2}  Let $V \in \rep L$, $\dim V = 2$. Then
	$V$ is indecomposable if and only if 	
	it is isomorphic to one, and only one, of $\U_{a,b}$ or to $\Vc_{a}^2$ or to $\W_{a}$.\end{prop}

\pf   
As discussed above, we may assume that
$V$ has a basis $\base$  such that  
$[g]_{\base} = \begin{pmatrix} a & 0 \\ 0 & -a\end{pmatrix}$ and
$[y]_{\base} = \begin{pmatrix} 0 & b \\ c & 0 \end{pmatrix}$ for some $b,c \in \ku$.
Now \eqref{eqrep} holds. We have four cases:
$b=c=0$: $V \simeq \ku_{a} \oplus \ku_{-a}$;  $c=0\neq b$: $V\simeq \W_{a}$;  $b=0 \neq c$: $V\simeq \W_{-a}$; 
$b ,c \in \ku^{\times}$:  $V \simeq \U_{a,bc}$. \epf

\subsubsection{Dimension 3}
Assume that $\dim V = 3$ and let $a\in \ku^{\times}$. 

\begin{prop}\label{prop:L-dim3} Let $V \in \rep L$, $\dim V = 3$. Then
$V$ is indecomposable if and only if 
$V$ is isomorphic to one, and only one, of $\Vc_{a}^{3}$, $\cC^3_{a}$,  $\Ec^3_{1,a}$, $\Ec^3_{2,a}$ 
or  $\Ec^3_{3,a,b}$, for unique $a,b \in \ku^{\times}$. \end{prop}

\pf  
If the action of $g$ is semisimple, then $V \simeq \cC^3_{a}$
by Lemmas \ref{lema:g-semisimple-1} and \ref{lema:g-semisimple-2}. 
Otherwise  Lemma \ref{lema:two-blocks-p=1} applies, since $V$  should have  a basis 
$\base$ where the action of $g$ is  given by
$[g]_{\base} ={\small \begin{pmatrix}
\jordan_{2}(a) & 0\\
0 & -a \end{pmatrix}}$.
\epf 

\subsubsection{Dimension 4}
Assume that $\dim V = 4$ and let $a\in \ku^{\times}$. 

\begin{prop}\label{prop:L-dim4} Let $V \in \rep L$, $\dim V = 4$. Then
	$V$ is indecomposable if and only if 
	$V$ is isomorphic to one, and only one, of $\Vc_{a}^4$, $\cD^4_{1,a}$, $\cD^4_{2,a,b}$, $\cD^4_{3,a}$, $\Ec^4_{1,a}$, $\Ec^4_{2,a}$, $\Ec^4_{3,a,b}$, $\Fc^{4}_{1,a}$, $\Fc^{4}_{2,a}$, $\Fc^{4}_{5,a,b,c}$,  $\Fc^{4,1}_{6,a,b,c}$, $\Fc^{4,1}_{8,a,b,c}$, $\cG^{4}_{a}$, $\cH^{4}_{1,a}$, $\cH^{4}_{2,a}$ or $\cH^{4}_{3,a,b}$ for unique $a,b,c \in \ku^{\times}$. \end{prop}
\pf  

\begin{align*}
\begin{tabular}{c | c|  c | c }
\noalign{\smallskip}    $A$  &  $B$ & Indecomposable & Reference
\\ \hline
$a\Id_3$   & $-a$  & none  &   Lemma \ref{lema:g-semisimple-2}
\\ \hline
$a\Id_2$   & $-a\Id_2$  & $\cD^{4}_{1,a}$, $\cD^{4}_{2,a,b}$ or  $\cD^{4}_{2,a,b}$  &   Lemma \ref{lema:g-semisimple-3}
\\ \hline  
$\jordan_{3}(a)$   & $-a$  & $\Ec^4_{1,a}$, $\Ec^4_{2,a}$ 
or  $\Ec^4_{3,a,b}$  &   Lemma \ref{lema:two-blocks-p=1}
\\\hline
$\jordan_{2}(a)$   & $\jordan_{2}(-a)$  &  $\Fc^{4}_{1,a}$, $\Fc^{4}_{2,a}$, $\Fc^{4}_{5,a,b,c}$,  $\Fc^{4,1}_{6,a,b,c}$ or $\Fc^{4,1}_{8,a,b,c}$  & Lemma \ref{lema:two-blocks-p=2}   
\\ \hline
$\jordan_{2}(a)$   & $-a\Id_2$  & $\cG^4_{a}$  &   Lemma \ref{lema:A-block-B-id-p>1}
\\ \hline
${\small \begin{pmatrix}
\jordan_{2}(a) & \\ & a
	\end{pmatrix}}$   & $-a$  & $\cH^{4}_{1,a}$, $\cH^{4}_{2,a}$ or $\cH^{4}_{3,a,b}$  &  Lemma \ref{lema:jordan-plus-1dim}
\\ \hline
\end{tabular}
\end{align*}
\epf 

\subsubsection{Dimension 5}
Assume that $\dim V = 5$ and let $a\in \ku^{\times}$. 

\begin{prop}\label{prop:L-dim5} Let $V \in \rep L$, $\dim V = 5$. Then
	$V$ is indecomposable if and only if 
	$V$ is isomorphic to one, and only one, of $\Vc_{a}^5$, $\cD^{5}_{4,a}$,   $\Ec^5_{1,a}$, $\Ec^5_{2,a}$, $\Ec^5_{3,a,b,c}$, $\Fc^{5}_{j,a}$, $j\in \I_{3}$, $\Fc^{5,i}_{4,a}$, $\Fc^{5}_{j,a,b,c}$, $j=5,7$, $\Fc^{5,i}_{j,a,b,c}$, $j=6,8$, $\cG^5_{a}$, $\cH^{5}_{1,a}$, $\cH^{5}_{2,a}$, $\cH^{5}_{3,a,b}$ or $\mathcal{I}^5_{a}$ for unique $a,b,c \in \ku^{\times}$. \end{prop}

\pf  
\begin{align*}
\begin{tabular}{c | c|  c | c }
\noalign{\smallskip}    $A$  &  $B$ & Indecomposable & Reference
\\ \hline
$a\Id_4$   & $-a$  & none  &   Lemma \ref{lema:g-semisimple-2}
\\ \hline
$a\Id_3$   & $-a\Id_2$  & $\cD^{5}_{4,a}$    &   Lemma \ref{lema:g-semisimple-3}
\\ \hline 
$\jordan_{4}(a)$   & $-a$  & $\Ec^5_{1,a}$, $\Ec^5_{2,a}$ 
or  $\Ec^5_{3,a,b}$  &   Lemma \ref{lema:two-blocks-p=1}
\\ \hline
$\jordan_{3}(a)$   & $\jordan_{2}(-a)$  & $\begin{array}{c}\Fc^{5}_{j,a},\, j\in \I_{3},\, \Fc^{5,i}_{4,a},\, \Fc^{5}_{j,a,b,c}, \\[.3em]
j=5,7,  \text{ or } \Fc^{5,i}_{j,a,b,c},\, j=6,8 \end{array}$&  Lemma \ref{lema:two-blocks-p=2}     
\\ \hline
$\jordan_{3}(a)$   & $-a\Id_2$  & $\cG^5_{a}$  &   Lemma \ref{lema:A-block-B-id-p>1}
\\ \hline
${\small \begin{pmatrix}
	\jordan_{3}(a) & \\ & a
	\end{pmatrix}}$   & $-a$  & $\cH^{5}_{1,a}$, $\cH^{5}_{2,a}$ or $\cH^{5}_{3,a,b}$ &  Lemma \ref{lema:jordan-plus-1dim}
\\ \hline
${\small \begin{pmatrix}
	\jordan_{2}(a) & \\ & \jordan_{2}(a)
	\end{pmatrix}}$   & $-a$  &   $\mathcal{I}^5_{a}$  &   Lemma \ref{lema:A-t-bloques-r-B-id-p1}
\\ \hline
\end{tabular}
\end{align*}
\epf

\subsection{Tensor products of some indecomposable modules}\label{subsec:tensor-simples-L}
Here we start the study of the tensor category $\rep L$. Let $\mathfrak{F}$ be a family of isomorphism classes of
indecomposable modules. Let  $\rep_{\mathfrak{F}} L$  be the full subcategory of $\rep L$
whose objects are direct sums of representatives of $\mathfrak{F}$, or in other words whose indecomposable components belong to $\mathfrak{F}$.

\begin{prop}
Let $\mathfrak{F}$ be the family of classes of the modules $\ku_{a}$, $\U_{a,b}$, $\W_{a}$, 
$a,b\in \ku^{\times}$.  Then $\rep_{\mathfrak{F}} L$ is a monoidal subcategory of $\rep L$.
\end{prop}

\pf Let $a,b,c,d\in \ku^{\times}$.
First we claim the results summarized in the following table;
the proofs are either straightforward  or else appear below.

\begin{align}\label{table:tensor}
\begin{tabular}{c |  c c c}
\noalign{\smallskip}  &  $\ku_{b}$  &  $\mathcal{U}_{1,d}$ & $\W_{1}$
\\ \hline
$(\underline{\;\;})^*$   & $\ku_{b^{-1}}$  & $\U_{1,-d}$ & $\W_{-1}$  
\\ \hline
$\ku_{a} \otimes \underline{\;\;}$ &  $\ku_{ab}$ & $\U_{a,a^{2}d}$  & $\W_{a}$  
\\ \hline
$\mathcal{U}_{1,c} \otimes \underline{\;\;}$  & $\mathcal{U}_{b,c}$  & $\U_{1,c+ d}\oplus \U_{1, c+ d}$, if $c \neq -d$ & 
$\U_{1,c}\oplus \U_{1,c}$
\\ 
& & $\W_{1}\oplus \W_{-1}$, if $c = -d$ & 
\\ \hline
$\W_{1} \otimes \underline{\;\;}$   & $\W_{b}$  & $\U_{-1,d}\oplus \U_{-1,d}$ &   $\W_{1} \oplus \W_{-1}$
\end{tabular}
\end{align}
Next, we claim that \eqref{table:tensor} implies the Proposition. For instance,
\begin{align*}
&\mathcal{U}_{a,c} \otimes \mathcal{U}_{b,d} \simeq \mathcal{U}_{1,c} \otimes \ku_a \otimes \mathcal{U}_{1,d} \otimes \ku_b
\simeq  \mathcal{U}_{1,c} \otimes \U_{a,a^{2}d} \otimes \ku_{b}
\\ &\simeq  \mathcal{U}_{1,c} \otimes \U_{1,a^{2}d} \otimes \ku_{ab} \simeq 
\begin{cases}
\U_{ab,c+ a^{2}d}\oplus \U_{ab, c+ a^{2}d},& c \neq a^{2}d,
\\		
\W_{-ab}\oplus \W_{ab},& c = a^{2}d.
\end{cases}
\end{align*}
\epf

We notice that $\rep_{\mathfrak{F}} L$ has interesting monoidal subcategories by considering some subsets of parameters,
e.~g. $a,b,c,d$ in a subring of $\ku$, with $a,b$ units. Also notice that $\rep_{\mathfrak{F}} L = \oplus_{\lambda\in \ku^{\times}} \rep_{\mathfrak{F}}^{\lambda} L$, where
$\rep^{\lambda}_{\mathfrak{F}} L := \rep_{\mathfrak{F}} L \cap \rep^{\lambda} L$. In particular, $\rep^{1}_{\mathfrak{F}} L$ appears to be interesting.

\begin{obs}
The monoidal category $\rep_{\mathfrak{F}} L$ can be interpreted as follows.
\begin{enumerate}[leftmargin=*,label=\rm{(\roman*)}]
	\item\label{item:monsubcat1}  Let $\hopfa$ be a Hopf algebra with the Chevalley property, i.e. the tensor product of any two simple $\hopfa$-modules is semisimple. Then the full subcategory $\reps \hopfa$ of $\rep \hopfa$ consisting of semisimple modules is monoidal.
Clearly, is the category of comodules over the Hopf subalgebra of the Sweedler dual of $\hopfa$ generated by the matrix coefficients of simple modules. But it is not a Serre subcategory of $\rep \hopfa$, in general.
For instance, if $\Gamma$ is an abelian group and $\hopfa = \ku\Gamma$, then $\reps \hopfa$ is   the category of
$\ku\widehat{\Gamma}$-comodules.

\item\label{item:monsubcat2} Now assume that $\hopfa$ is a Hopf subalgebra of a Hopf algebra $\hopfb$ 
and let $\mathfrak C$ be a monoidal subcategory of $\rep \hopfa$. Then the full subcategory 
$\rep_{\mathfrak C} \hopfa$ of $\rep \hopfb$ consisting of $\hopfb$-modules that when 
restricted to $\hopfa$ belong to $\mathfrak{C}$ is monoidal.
\end{enumerate}

Then $\rep_{\mathfrak{F}} L$ is a monoidal subcategory of $\rep_{\mathfrak C} L$, where 
$\mathfrak C$ is $\reps \ku G$.
\end{obs}

\subsubsection{Two simple modules, $\dim 2$} 
Let $c,d \in \ku^{\times}$ and $V=\U_{1,c}\otimes \U_{1,d}$.
Let $v_1,v_2$ be a  basis of $\U_{1,c}$, $w_1,w_2$ a basis of $\U_{1,d}$, both realizing \eqref{simple:2dim}. 
In the basis $u_1=v_1\otimes w_1$, $u_2=v_2\otimes w_2$, $u_3=v_1\otimes w_2$, $u_4=v_2\otimes w_1$, the action is
\begin{align}\label{eq:simple-tensor-product-action}
\begin{aligned}
&gu_1=u_1,&  &yu_1=u_3+u_4,& &gu_3=-u_3,& &yu_3= du_1+u_2, \\
&gu_2=u_2,&  &yu_2=c u_3-d u_4,& &gu_4=-u_4,& &yu_4= cu_1-u_2.
\end{aligned}
\end{align}
We claim that
\begin{align}\label{eq:simple-tensor-product} 
	\begin{aligned}
c+ d &\neq 0 \implies
\U_{1,c}\otimes \U_{1,d}\simeq \U_{1,c+ d}\oplus \U_{1, c+ d},
\\		
c+ d &= 0 \implies
\U_{1,c}\otimes \U_{1,d}\simeq \W_{-1}\oplus \W_{1}.
	\end{aligned}
\end{align}
\pf $c+ d \neq 0$: Here 
$u_{1}$, $u_3+u_4$, $u_{2}$, $c u_3-d u_4$ is a basis of $V$. By
\eqref{eq:simple-tensor-product-action}, $\langle u_{1}, u_3+u_4 \rangle \leq V$ and 
$\langle u_{2}, c u_3-d u_4 \rangle \leq V$, both  isomorphic to 
$\U_{1,c+d}$.\smallbreak

 $c+ d = 0$: Now $u_3+u_4$, $u_{1}$,  $d u_{1} + u_{2}$, $u_3$ is a basis of $V$. By
 \eqref{eq:simple-tensor-product-action}, $\W_{-1} \simeq \langle u_3+u_4, u_{1} \rangle \leq V$ and 
 $\W_{1} \simeq \langle d u_{1} + u_{2}, u_{3} \rangle \leq V$.\epf

\subsubsection{Simple and indecomposable, $\dim 2$} 
Let $c,d \in \ku^{\times}$ and $V=\U_{1,c}\otimes \W_{1}$.
Let $v_1,v_2$ be a  basis of $\U_{1,c}$, $w_1,w_2$ a basis of $\W_{1}$. 
In the basis $u_1=v_1\otimes w_1$, $u_2=v_2\otimes w_2$, $u_3=v_1\otimes w_2$, $u_4=v_2\otimes w_1$
of $V$ the action is
\begin{align}\label{eq:simple-indec-tensor-product-action}
	\begin{aligned}
		&gu_1=u_1,&  &yu_1=u_4,& &gu_3=-u_3,& &yu_3= u_1+u_2, \\
		&gu_2=u_2,&  &yu_2=c u_3- u_4,& &gu_4=-u_4,& &yu_4= cu_1.
	\end{aligned}
\end{align}
We claim that
\begin{align}\label{eq:simple-indec-tensor-product} 
\U_{1,c}\otimes \W_{1} &\simeq \U_{1,c}\oplus \U_{1,c},&
\W_{1}\otimes \U_{1,d} &\simeq \U_{-1,d}\oplus \U_{-1,d}.
\end{align}

\noindent \emph{Proof.}  By \eqref{eq:simple-indec-tensor-product-action}, $\langle u_1, u_{4} \rangle \leq V$ and 
$\langle u_{2}, c u_3- u_4 \rangle \leq V$, both  $\simeq \U_{1,c}$, hence the first isomorphism. 
Next compute $\ku_{-1} \otimes(\U_{1,c}\otimes \W_{1})^* \simeq \ku_{-1} \otimes (\U_{1,c}\oplus \U_{1,c})^*$:
\begin{align*}
\ku_{-1} \otimes(\W_{1}^* \otimes \U_{1,c}^*) &\simeq \ku_{-1} \otimes( \W_{-1} \otimes  \U_{1,-c})
\simeq \ku_{-1} \otimes (\U_{1,-c}\oplus \U_{1,-c}) 
\\ & \implies 
 \W_{1} \otimes  \U_{1,-c} \simeq \U_{-1,-c}\oplus \U_{-1,-c}. \qquad \qed
\end{align*}

\subsubsection{Two indecomposable modules of dimension $2$}
Let $v_1,v_2$ and $w_1,w_2$ be basis of two copies of $\W_{1}$.
In the basis $u_1=v_1\otimes w_1$, $u_2=v_2\otimes w_2$, $u_3=v_1\otimes w_2$, $u_4=v_2\otimes w_1$
of $V=\W_{1}\otimes \W_{1}$ the action is
\begin{align}\label{eq:indec-indec-tensor-product-action}
\begin{aligned}
&gu_1=u_1,&  &yu_1=0,& &gu_3=-u_3,& &yu_3= u_1, \\
&gu_2=u_2,&  &yu_2= u_3- u_4,& &gu_4=-u_4,& &yu_4= u_1.
\end{aligned}
\end{align}
We claim that
\begin{align}\label{eq:indec-indec-tensor-product} 
\W_{1}\otimes \W_{1} &\simeq \W_{1} \oplus \W_{-1}.
\end{align}
\pf  By \eqref{eq:indec-indec-tensor-product-action}, $\W_{1} \simeq \langle u_1, u_{3} \rangle \leq V$ and 
$\W_{-1} \simeq \langle u_3- u_4, u_{2} \rangle \leq V$.
\epf

\section{The bosonization of the Jordan plane}\label{sec:bos-jordan}
Let $\Ac = \ku \langle \letra_1,\letra_2 \rangle$ modulo the ideal generated by the quadratic relation
\begin{align}\label{rel:quad}
	\letra_2\letra_1&-\letra_1\letra_2+\frac{1}{2}\letra_1^{2}.
\end{align}
This is the well-known \emph{Jordan plane}. Let $\Gletra$ be an infinite cyclic group denoted multiplicatively with a fixed generator $\letrag$.
Let  $\Vc= \ku \letra_1 \oplus \ku \letra_2\in \ydG$ with grading $\Vc_\letrag = \Vc$ and action $\letrag \cdot \letra_1 = \letra_1$, $\letrag \cdot \letra_2 = \letra_1 + \letra_2$.
Then $\Ac \simeq \toba(\Vc)$, cf. \cite[ Prop. 3.4]{aah-triang}.
Let $H:=\Ac\#\ku \Gletra$ the bosonization of $\Ac$ by $\ku \Gletra$; i.~e. $H = \ku \langle \letra_1,\letra_2, \letrag^{\pm 1} \rangle$ modulo the ideal generated by   \eqref{rel:quad}, 
$\letrag^{\pm}\letrag^{\mp} - 1$,
\begin{align}\label{rel:boso1}
	\letrag\letra_1\letrag^{-1}&-\letra_1,
	\\\label{rel:boso2}
	\letrag\letra_2\letrag^{-1} &- (\letra_1+\letra_2).
\end{align}
This is a Hopf algebra with the comultiplication determined by
\begin{align*}
	\Delta(\letrag^{\pm 1})=\letrag^{\pm 1}\otimes \letrag^{\pm 1}, & & \Delta(\letra_i)=\letra_i\otimes 1+\letrag\otimes \letra_i,& &i=1,2.&
\end{align*}
The set $\{\letra_1^a \letra_2^b \letrag^c:  a,b\in\N_0, c \in \mathbb Z\}$ is a basis of $H$, whose $\GK$ is 3.

\subsection{The Hopf algebra $\overline{H}$}\label{subsec:H-barra}
From the defining relations, we see that the left ideal $H \letra_1$ is a two-sided, as well as a Hopf ideal.
Then $\overline{H} := H / H \letra_1$ is the commutative Hopf algebra $\ku\la\overline{\letrag}^{\pm 1},\overline{\letra}_2\ra$, that is
\begin{align*}
	\overline{H} &  \simeq \Oc (\mathbb{B}),&
	\mathbb{B} &= \left\{{\small \begin{pmatrix} a & b \\ 0 & a^{-1} \end{pmatrix}}:(a,b) \in \ku^{\times} \times \ku \right\}
	\leq SL_2(\ku), 
\end{align*}
where $\Oc$ stands for the algebra of regular functions. 
Therefore, the tensor category $\rep \overline{H}$ reflects the group structure of $\mathbb{B}$.
The $\overline{H}$-modules, are described as follows. 
Let  $A\in \Gln_n(\ku)$ and $B\in \End(\ku^n)$ such that $AB=BA$. 
We denote by $\ku^n_{A,B} \in  \rep \overline{H}$ the vector space $\ku^n$ 
with action given by $\overline{\letrag} \mapsto A$ and $\overline{\letra}_2 \mapsto B$. 
Every $V\in  \rep \overline{H}$ is isomorphic to some $\ku^n_{A,B}$; also,
$\ku^n_{A,B}\simeq \ku^n_{A',B'}$ iff $A$, $B$ and $A'$, $B'$ are simultaneously conjugated.
For $n =1$, let $\ku_{\gamma} = \ku_{a,b}$, where $\gamma = {\small \begin{pmatrix} a & b \\ 0 & a^{-1} \end{pmatrix}} \in \mathbb{B}$; this says that the simple $\overline{H}$-modules are classified by the points of $\mathbb{B}$.
Given $a\in \ku^{\times}$, $b,c\in \ku$, we have $\overline{H}$-modules of dimension 2 given by
\begin{align*}
	\mathcal{J}_{a,b}:  \overline{\letrag} &\mapsto {\small \begin{pmatrix}
			a & 0 \\
			0 & a \end{pmatrix}},& \overline{\letra}_2& \mapsto {\small \begin{pmatrix}
			b & 1 \\
			0 & b \end{pmatrix}},&  
	\mathcal{K}_{a,b,c}: \overline{\letrag}& \mapsto {\small \begin{pmatrix}
			a & 1 \\
			0 & a \end{pmatrix}},&  \overline{\letra}_2& \mapsto {\small \begin{pmatrix}
			b & c \\
			0 & b \end{pmatrix}}.
\end{align*}
Clearly these are pairwise non-isomorphic indecomposable $\overline{H}$-modules.
We leave to the reader the (elementary) proof of the following result.

\begin{lema}\label{rep_dime_two}
	\begin{enumerate}[leftmargin=*,label=\rm{(\roman*)}]
		\item 
		If $V \in  \rep \overline{H}$ is indecomposable of dimension 2,
		then either $V \simeq\mathcal{J}_{a,b}$, or $V \simeq\mathcal{K}_{a,b,c}$ for unique $a,b,c$.
		\item $ \dim \ext_{\overline{H}}^1(\ku_\gamma,\ku_\eta) = \delta_{\gamma, \eta}$.
		\item \label{generatorH} 
		If $a \in \ku^{\times}$ and $b,c \in \ku$, then
		\begin{align*}
			\ku_{a,b-a} \otimes \mathcal{J}_{1,1} &\simeq \mathcal{J}_{a,b} \simeq \mathcal{J}_{1,1} \otimes \ku_{a,b-1},&
			\mathcal{J}_{a,b}^{*}&\simeq \mathcal{J}_{a^{-1},-ba^{-1}};
			\\
			\ku_{a,b-a} \otimes \mathcal{K}_{1,1,c}&\simeq \mathcal{K}_{a,b,c} \simeq \mathcal{K}_{1,1,ac-(b-1)} \otimes\ku_{a,b-1},&
			\mathcal{K}_{a,b,c}^{*} &\simeq  
			\mathcal{K}_{a^{-1},-ba^{-1},ca-b}. 
		\end{align*}
	\end{enumerate}\qed
\end{lema}

\subsection{The category $\rep H$}\label{subsec:simple}
The projection $H \to \overline{H}$ induces a functor \newline $\rep \overline{H} \to \rep H$. We carry over the notation 
along this functor.
Conversely, let $V \in \rep H$ and $V_{0}:=\ker \letra_1$, giving a functor $\rep H \to \rep \overline{H}$.
If $V = V_{0}$ has dimension $n$,   then $V  \in \rep \overline{H}$, hence $V \simeq \ku^n_{A,B}$
for some $A$, $B$.

\medbreak
Let $V \in \rep H$. Then $V_0 \leq V$ since $\letra_1 (\letra_2 V_0)=0$ by \eqref{rel:quad}, and $\letra_1(\letrag V_0)=0$ by  \eqref{rel:boso1}.

\begin{lema}\label{lema:iyudu} \cite[Lemma 2.1]{Na}
	If $U \in \rep \Ac$, then $\letra_1$ acts nilpotently on $U$; in particular 
	$V_{0} \neq 0$.  \qed
\end{lema}

The Lemma applies to $V \in \rep H$ via  the evident restriction functor.

\begin{prop} \label{prop:simple-H}
If $V \in \rep H$ is irreducible, then  $V \simeq\ku_{\gamma}$, for a unique  $\gamma  \in \mathbb{B}$. \qed
\end{prop}

We next classify the indecomposable $H$-modules of dimension $2$. 

\begin{prop}\label{prop:pol_min} Let $V \in \rep H$, $\dim V=2$. Then $V = V_0$.
	Therefore, if $V$ is indecomposable,
	then either $V \simeq\mathcal{J}_{a,b}$, or $V \simeq\mathcal{K}_{a,b,c}$ for unique $a,b,c$.
	Also
	\begin{align*}
		\dim \ext_{H}^1(\ku_\gamma,\ku_\eta) = \delta_{\gamma, \eta}.
	\end{align*}
\end{prop}
\pf If $V\neq V_0$ then there exists a basis $\Omega$ of $V$ such that
\begin{align*}
	& [\letra_1]_{\Omega}=\begin{pmatrix}
		0 & 1 \\
		0 &  0 \end{pmatrix},& [\letrag]_{\Omega}=\begin{pmatrix}
		a & b \\
		c &  d \end{pmatrix}, & & [\letra_2]_{\Omega}=\begin{pmatrix}
		e & f \\
		h &  i \end{pmatrix}.&
\end{align*}
By \eqref{rel:boso1},  $c=0$ and $a=d$. From $\letra_1\letra_2=\letra_2\letra_1$ we see that $h=0$ and $e=i$. By \eqref{rel:boso2}, we conclude  that $af+be=a(f+1)+be$. Hence, $a=0$ and $\letrag$ is not invertible,  a contradiction.
The last claims follow from Lemma \ref{rep_dime_two}. \epf

\section{The bosonization of the super Jordan plane}\label{sec:superjordan}
Let $x_{21}=x_1x_2+x_2x_1$ in the free associative algebra in generators $x_1$ and $x_2$. Let $\mathcal{B}$ be the algebra generated by $x_1$ and $x_2$ with defining relations
\begin{align}\label{eq:rels-B(V(-1,2))-1}
&x_1^2, \\
\label{eq:rels-B(V(-1,2))-2}
&x_2x_{21}- x_{21}x_2 - x_1x_{21}.
\end{align}
The algebra $\mathcal{B}$, introduced in \cite{aah-triang}, is called the \emph{super Jordan plane}. 
Let $\Vc' = \ku x_1 \oplus \ku x_2\in \ydg$ with grading $\Vc'_g = \Vc'$ and action $g \cdot x_1 = -x_1$, $g \cdot x_2 = -x_1 + x_2$.
Then $\mathcal{B} \simeq \toba(\Vc')$, cf. \cite[ Prop. 3.5]{aah-triang}.
Let $K:=\mathcal{B}\#\ku G$ the bosonization of $\mathcal{B}$ by $\ku G$; i.~e. $K = \ku \langle x_1,x_2, g^{\pm 1} \rangle$ modulo the ideal generated by  \eqref{eq:rels-B(V(-1,2))-1}, \eqref{eq:rels-B(V(-1,2))-2}, $g^{\pm}g^{\mp} - 1$,
\begin{align}\label{rel:superboso1}
gx_1g^{-1}& + x_1,
\\\label{rel:superboso2}
gx_2g^{-1} &- x_1+x_2.
\end{align}
This is a Hopf algebra with the comultiplication determined by
\begin{align*}
\Delta(g^{\pm 1})=g^{\pm 1}\otimes g^{\pm 1}, & & \Delta(x_i)=x_i\otimes 1+g\otimes x_i,& &i=1,2.&
\end{align*}
The set $\{x_1^a x_{21}^bx_2^cg^d:   a\in \I_{0,1}, b, c\in\N_0, d \in \mathbb Z\}$ is a basis of $K$, whose $\GK$ is 3.
The following identities are valid in $K$:
\begin{align}\label{eq:relations1}
x_{21}x_1&=x_1x_{21},
\\\label{eq:relations2} x_2^2x_1 &=x_1x_2^2+x_1x_2x_1,
\\\label{eq:relations3} x_{21}x_2^2 &=(x_2^2-x_{21})x_{21},
\\\label{eq:relations4} gx_{21} &= x_{21}g,
\\\label{eq:relations5} gx_{2}^2 &=(x_{2}^2 - x_{21})g.
\end{align}

\begin{lema}\label{lem:hopfideal} 
\begin{enumerate}[leftmargin=*,label=\rm{(\roman*)}]
\item\label{item:hopfideal1}  There is an injective algebra map $\varphi:H \to K$ given by
\begin{align}\label{eq:Jordan-to-superJordan}
\varphi(\letra_1)&= x_{21},& \varphi(\letra_2)&= -\frac{1}{2}x_{2}^2,& \varphi(\letrag)&= g^2.
\end{align}
\item\label{item:hopfideal2} $Kx_1K=Kx_1+Kx_{21}$ is a Hopf ideal and $ K/Kx_1K \simeq L$.
\end{enumerate}
\end{lema}

\pf \ref{item:hopfideal1}: It is not difficult to check that $\varphi$ is well-defined; indeed
\eqref{rel:quad} follows from  \eqref{eq:relations3}, \eqref{rel:boso1} from  \eqref{eq:relations4} and \eqref{rel:boso2}  from  \eqref{eq:relations5}. The injectivity is verified using the PBW-bases.

\ref{item:hopfideal2}: Let $r=x_2-x_1$. We prove by induction on $n$  that 
\begin{align}\label{eq:auxiliar1-superjordan}
x_1x_2^{n}= (-1)^nx_2^nx_1+ \left((-1)^{n-1}x_2^{n-1}+ \cdots + x_2^2r^{n-3}-x_1r^{n-2}\right) x_{21}.
\end{align}
For $n=1$, \eqref{eq:auxiliar1-superjordan} is trivial. 
If $n >1$ and  \eqref{eq:auxiliar1-superjordan} holds for $n-1$, then $x_1x_2^{n} =$
\begin{align*}
&=\left\{(-1)^{n-1}x_2^{n-1}x_1+\left((-1)^{n-2}x_2^{n-2}+ \cdots + x_2^2r^{n-4}-x_1r^{n-3}\right)x_{21}\right\}x_2\\
&= (-1)^{n-1}x_2^{n-1}x_1x_2+\left((-1)^{n-2}x_2^{n-2}+ \cdots + x_2^2r^{n-4}-x_1r^{n-3}\right)x_{21}x_2\\
&= (-1)^{n-1}x_2^{n-1}(x_{21\!\! }-\! x_2x_1)\!\! +\! \left((-1)^{n-2}x_2^{n-2}+\! \cdots\! + x_2^2r^{n-4}-x_1r^{n-3}\right)rx_{21}\\
&= (-1)^nx_2^nx_1+\left((-1)^{n-1}x_2^{n-1}+ \cdots + x_2^2r^{n-3}-x_1r^{n-2}\right)x_{21}.\end{align*}
Hence $x_1x_2^{n} \in Kx_1+Kx_{21}$, for all $n \in \mathbb{N}$, and $x_1K \subset Kx_1+K_{21}$.  
The isomorphism is verified using the PBW-bases.
\epf

\subsection{Relations between $\rep K$, $\rep L$ and $\rep H$} Let $V \in \rep K$.
If $x_1 = 0$ in $V$, then also $x_{21} =0$; by Lemma \ref{lem:hopfideal}, we conclude that  $V \in \rep L$,
a category discussed in \S \ref{sec:borel}. Hence we may assume that $x_1 \neq 0$.

\smallbreak
More generally, since $x_1^2 = 0$, we may think on it as a differential  on $V$ and consider its homology. 
Namely, define
\begin{align}\label{eq:functors}
\Kg(V) &= \ker x_1 \text{ on }V,& \Igo(V) &=   x_1(V),& \Hg(V) &= \Ig(V)/ \Kg(V).
\end{align}

\begin{prop}\begin{enumerate}[leftmargin=*,label=\rm{(\roman*)}]
\item\label{item:functors1}   $\Kg, \Igo, \Hg$ are $\ku$-linear functors $\rep K \to \rep H$.

\item\label{item:functors2} $\Kg$ is left exact.

\item\label{item:functors3} $\Igo$ is right exact.
	\end{enumerate}
\end{prop}

\pf \ref{item:functors1}: The claims for $\Kg, \Igo$ are consequences of \eqref{eq:relations1} and \eqref{eq:relations2},
and imply in turn that of $\Hg$.   
\ref{item:functors2} and \ref{item:functors3} are standard.   \epf

\subsection{Simple modules}\label{subsec:simple_super}
We show that the  classification of the simple objects in  $\rep K$ reduces to those in $\rep L$
given in Proposition \ref{prop:simple_boso_super}.

\begin{theorem}\label{theor:simple_boso_super}
Let $V\in \rep K$ irreducible. Then $V\in \rep L$, in particular  $\dim V=1$ or $\dim V=2$.
\end{theorem}

\pf First we claim that $W:=\ker x_{21}$ is a submodule of $V$.
By \eqref{eq:relations2} and \eqref{eq:relations4}, $W$ is stable by $x_1$ and $g$. Let $u\in W$ so that
$x_1x_2u= -x_2x_1u$.  Then 
\begin{align*}
&x_{21}x_2u=x_1x_2^2u+x_2 x_1x_2u=x_1x_2^2u-x_2^2x_1u 
 \overset{\eqref{eq:relations2}}{=}-x_1x_2x_1u= x_1^2 x_2 u = 0,
\end{align*}
hence $ x_2u\in W$. By Lemmas \ref{lema:iyudu} and \ref{lem:hopfideal} \ref{item:hopfideal1}, $W \neq 0$ 
hence $W = V$, i.~e. $x_1x_2=-x_2x_1$ in $V$. Thus $0\neq \ker x_1$ is $x_2$-stable and by \eqref{rel:superboso1}, $\ker x_1 \leq V$; consequently $\ker x_1=V$. 
\epf

\subsection{Indecomposable modules} For convenience, we set
\begin{align*}
s &=x_{21}, & t&=x_2^2.
\end{align*}
 Let $V \in \rep K$ such that $x_1 \neq 0$ on $V$. 
 By Lemmas \ref{lema:iyudu} and \ref{lem:hopfideal}, $s$ is nilpotent.
The aim of this Subsection is to establish the following result.

\begin{prop}\label{prop:eigen}  Let $V \in \rep K$ such that $x_1 \neq 0$. 
	
\begin{enumerate}[leftmargin=*,label=\rm{(\roman*)}]
\item\label{it:SJ-eigent1}  $V_t^{(a)} \leq V$, $a\in \ku$; hence $V = \oplus_{a\in \ku^{\times}} V_t^{(a)}$ is a decomposition in $\rep K$.

\item\label{it:SJ-eigent2} If $V$ is indecomposable, then $V =V_t^{(a)}$, 
for a unique eigenvalue  $a$ of $t$. Hence $\spec x_2 \subseteq \{\pm \nu\}$, where $\nu^2 = a$.

\item\label{it:SJ-eigeng1} Let $V[\lambda]:=V_g^{(b)}\oplus V_g^{(-b)}$, $\lambda = b^2\in \ku^{\times}$.
Then $V[\lambda] \leq V$, hence  $V = \oplus_{\lambda\in \ku^{\times}} V[\lambda]$ is a decomposition in $\rep K$.

\item\label{it:SJ-eigeng2} If $V$ is indecomposable, then $V = V[\lambda]$, 
for a unique $\lambda = b^2\in \ku^{\times}$ and $\spec g = \{\pm b\}$.
\end{enumerate}	
\end{prop}

We start by relations in $\cB$ that might be useful for other problems. Set
\begin{align*}
\zeta_{n,j} &=\frac{n!}{(n-j)!},& j, n &\in \N_0,\ j\leq n.
\end{align*}

\begin{lema}\label{lem:other-rel} Let $a, b \in \ku$ and set $w =s-a,\, z =t-a \in \cB$.
Then for every $n \in \N$, we have
\begin{align}
\label{eq:clubsuit}
z^nx_1 &= x_1\sum_{j\in \I_{0,n}} \zeta_{n,j} s^jz^{n-j}
\\ \label{eq:circledast} 
z^{n}g&= g(z+s)^n, 
\\\label{eq:heartsuit}
(z+s)^n &=\sum_{j\in \I_{0,n}}\zeta_{n,j}s^jz^{n-j},
\\\label{eq:lozenge}
w^nx_2 &=x_2w^n-nx_1sw^{n-1}.&
\\ \label{eq:diamondsuit}
(g-b)^nx_1 &= (-1)^nx_1(g+b)^n, 
\\\label{eq:diamondsuitdiamondsuit}
(g-b)^{n}x_2&=(-1)^{n-1} \left(nx_1g-x_2(g+b) \right) (g+b)^{n-1}.
\end{align} 
	
\end{lema}
\noindent \emph{Proof.} 
\eqref{eq:clubsuit} is \cite[Lemma 2.9]{rep-super}.
\eqref{eq:circledast}: $x_2g=-g(x_1+x_2)$ by \eqref{rel:superboso1} and \eqref{rel:superboso2}. Hence $tg=x_2^2g=g(x_1+x_2)^2=g(x_{21}+x_2^2)=g(s+t)$. Thus $zg=g(z+s)$ and $z^ng=g(z+s)^n$. 

\smallbreak
\noindent\eqref{eq:heartsuit}: We proceed by induction on $n$; the case $n=1$ is trivial. Let $n>1$. By \eqref{eq:relations3}, $zs=sz+s^2$ whence $zs^n=s^nz+ns^{n+1}$. Then
\begin{align*}
(z+s)^n &=(z+s)\sum\limits_{j=0}^{n-1}\zeta_{n-1,j}s^jz^{n-1-j}
 =\sum\limits_{j=0}^{n-1}\zeta_{n-1,j}\left(zs^jz^{n-1-j} + s^{j+1}z^{n-1-j}\right)\\
& =\sum\limits_{j=0}^{n-1}\zeta_{n-1,j} \left(s^jz^{n-j}+ js^{j+1}z^{n-1-j}+ s^{j+1}z^{n-1-j}\right)\\
& =z^n+\sum\limits_{k=1}^{n-1}(\zeta_{n-1,k} +k\zeta_{n-1,k-1})s^kz^{n-k}+ n!s^n
=\sum\limits_{k=0}^{n}\zeta_{n,k}s^kz^{n-k}.
\end{align*}

\noindent\eqref{eq:lozenge}: also by induction on $n$. For $n=1$, we have 
\begin{align*}
wx_2&=x_1x_2^2+x_2x_1x_2-a x_2
\overset{\eqref{eq:relations2}}{=} x_2^2x_1-x_1x_2x_1+x_2x_1x_2-a x_2\\
&= x_2(s-a\id)-x_1s = x_2w-x_1s.
\end{align*}
Let $n>1$ and suppose that the identity is true for $n-1$. Then 
\begin{align*}
w^nx_2&=w(x_2w^{n-1}-(n-1)x_1sw^{n-2}) 
\\ &=(x_2w-x_1s)w^{n-1}-(n-1)x_1sw^{n-1}
=x_2w^n-nx_1sw^{n-1}. 
\end{align*}	

\noindent\eqref{eq:diamondsuit} follows at once from \eqref{rel:superboso1}. 
\eqref{eq:diamondsuitdiamondsuit}:
For $n=1$, we have  
\begin{align*}
(g-b)x_2=gx_2-b x_2\overset{\eqref{rel:superboso2}}{=}(x_1-x_2)g-b x_2=x_1g-x_2(g + b).
\end{align*}
For $n>1$ we compute 
\begin{align*}
(g-b)^nx_2 &= (-1)^{n-2}(g-b) \left\{ (n-1)x_1g-x_2(g+b) \right\} (g + b)^{n-2}\\
&=(-1)^{n-2} \left\{ -(n-1)x_1g -x_1g + x_2(g+b)  \right\} (g+b)^{n-1}\\
&=(-1)^{n-1} \left(nx_1g-x_2(g+b) \right) (g+b)^{n-1}. \qquad\qquad\qquad\qquad \qed
\end{align*} 

\noindent\emph{Proof of Proposition \ref{prop:eigen}.} \ref{it:SJ-eigent1}:
Clearly, $V_t^{(a)}$ is stable by $x_2$. By \eqref{eq:clubsuit}, \eqref{eq:circledast} and \eqref{eq:heartsuit},
 with $n = 2 \dim V$, it is stable by $x_1$ and $g$. \ref{it:SJ-eigent2} follows from \ref{it:SJ-eigent1}.

\ref{it:SJ-eigeng1}: By \eqref{eq:diamondsuit}, respectively \eqref{eq:diamondsuitdiamondsuit},
 $x_1\cdot V_g^{(\pm b)}\subset V_g^{(\mp b)}$ and $x_2\cdot V_g^{(\pm b)}\subset V_g^{(\mp b)}$. Hence $V[\lambda] \leq V$, implying the first claim in \ref{it:SJ-eigeng2}.
By \eqref{eq:diamondsuit}, since $x_1\neq 0$, there exists  $b \in \spec g$ such that $V_g^{(-b)}\neq 0$, thus $-b \in \spec g$.
\qed

\bigbreak  Let $V \in \rep K$ and $\lambda = a^2\in \ku^{\times}$ such that $V = V[\lambda]$.
Let $\base_{\pm}$ be a basis of $V_{g}^{(\pm a)}$ and $\base:=\base_+\cup \base_{-}$;
let $\ell = \dim V_{g}^{(a)}$, $\wp = \dim V - \ell$.  Then 
\begin{align} \label{rep:general}
[g]_{\base} &=\begin{pmatrix}
A & 0 \\
0 &  B \end{pmatrix}, & 
[x_1]_{\base} &= \begin{pmatrix}
0 & C  \\
D & 0 \end{pmatrix},&
[x_2]_{\base} &=\begin{pmatrix}
0 & E  \\
F & 0 \end{pmatrix},
\end{align}
where $A\in \GL_{\ell}(\ku)$, $B\in \GL_{\wp}(\ku)$, $C,E\in\operatorname{M}_{\ell\times \wp}(\ku)$ and $D,F\in\operatorname{M}_{\wp\times \ell}(\ku)$. 
Conversely,  $g$, $x_1$ and $x_2$ given by \eqref{rep:general} define a representation of $K$ 
if and only if
\begin{align} \label{matrices}
	CD&=0, & DC&=0, & CF(C+E)&=EFC, \\
	\label{matrices1} AC&=-CB,&  BD&=-DA,& (C-E)B&=AE,\\
	\label{matrices2} DE(D+F)&=FED,& (D-F)A&=BF. &&
\end{align}

Again we seek to describe conditions that guarantee that $V$ is indecomposable, 
assuming that $A$ and $B$ are in Jordan form. Actually this will be done in some 
special cases.

\begin{obs} Let  $V\in \rep K$ indecomposable given by \eqref{rep:general}, $\dim V =n \geq 2$. If $A=a\id_{\ell}$ and $B=-a \id_{\wp}$, then $x_1=0$. In fact, from $AE-(C-E)B=0$ follows $aC=0$ hence $C=0$. Similarly, $BF-(D-F)A=0$ implies $D=0$. Particularly, the indecomposable $K$-modules of dimension $2$ are just the indecomposable $L$-modules
by Proposition \ref{prop:eigen} \ref{it:SJ-eigeng2}.
\end{obs}

\subsubsection{Representations with $A$ and $B$ Jordan blocks}
By the same arguments used in Lemma \ref{lema:two-blocks} we have the following. 

\begin{lema}\label{lema:two-blocks-case-super} Let $V\in \rep K$ given by \eqref{rep:general} where
	$A$ and $B$ are Jordan blocks. Then $V$ is indecomposable if and only if $x_1 \neq 0$ or $x_2 \neq 0$. \qed
\end{lema}

Assume that $V$ is a $K$-module of dimension $n > 3$ such that $x_1\neq 0$. Define the representations 
$\mathcal{L}^n_{1,a,b}$ and $\mathcal{L}^n_{2,a,b}$, $a \in \ku^{\times}$, $b \in \ku$ by \eqref{rep:general} where $A=\jordan_{n-1}(a)$, $B=-a$ and $x_1$, $x_2$ acts as follows:
\begin{align*}
&\mathcal{L}^n_{1,a,b}: & &  
x_1 \mapsto \begin{pmatrix} 
	0 & 0\\
	e_{n-1}^{t} & 0 \\ \end{pmatrix}, & &
x_2 \mapsto \begin{pmatrix} 
	0 & be_1\\
	ae_{n-2}^t & 0 \end{pmatrix}. \\
&\mathcal{L}^n_{2,a,b}: & & 
x_1 \mapsto \begin{pmatrix} 
	0 & e_1\\
	0 & 0  \end{pmatrix}, & &
x_2 \mapsto \begin{pmatrix} 
	0 & -ae_2 \\
	be_{n-1}^t & 0  \end{pmatrix}. \end{align*}
Clearly these are pairwise non-isomorphic indecomposable $K$-modules.

\begin{lema}\label{lema:two-blocks-p=1-super-result} Let $n > 3$ and $V\in \rep K$  given by \eqref{rep:general} where
	$A =  \jordan_{n-1}(a)$ and $B = -a$. 
	Then the following are equivalent:
	\begin{enumerate}[leftmargin=*,label=\rm{(\roman*)}]
		\item\label{item:two-blocks-p=1-1-super} $V$ is indecomposable,
		\item\label{item:two-blocks-p=1-2-super} $x_1 \neq 0$ or $x_2\neq 0$,
		\item\label{item:two-blocks-p=1-3-super} $V$ is isomorphic to one, and only one, of $\mathcal{L}^{n}_{1,a,b}$ or $\mathcal{L}^{n}_{2,a,b}$, for unique $a \in \ku^{\times}$, $b \in \ku$.
	\end{enumerate}
\end{lema}
\pf 
 (i) $\Leftrightarrow$ (ii) is Lemma \ref{lema:two-blocks-case-super}. We prove that (ii) $\Leftrightarrow$ (iii). Here $C=\sum_{i=1}^{n-1}c_ie_i$, $D=\sum_{i=1}^{n-1}d_ie_i^t$, $E= \sum_{i=1}^{n-1}f_ie_i$ and $F= \sum_{i=1}^{n-1}h_ie_i^t$ with $c_i,d_i, f_i$, $h_i \in \ku$, $i \in \I_{n-1}$. By \eqref{matrices}, \eqref{matrices1} and \eqref{matrices2} we have two possibilities: 
 \begin{itemize}[leftmargin=*]\renewcommand{\labelitemi}{$\circ$}
 	\item  $C=0$, $D=de^t_{n-1}$, $E=fe_1$, $F=ade^t_{n-2}+he^t_{n-1}$, with $d \in \ku^{\times}$, $f,h\in \ku$. 
 	\item  $C=ce_{1}$, $D=0$, $E=fe_1-ace_2$, $F=he^t_{n-1}$, with $c \in \ku^{\times}$, $f,h\in \ku$. 
 \end{itemize}
In the first, the basis $e_1, -h(ad)^{-1}e_1+e_2, \dots, -h(ad)^{-1}e_{n-2}+e_{n-1}, de_n$ gives $V \simeq \mathcal{L}^{n}_{1,a,df} $. In the second, the basis $ce_1, -fa^{-1}e_1+ce_2, \dots,$ $ -fa^{-1}e_{n-2}+ce_{n-1}, e_n$  gives $V \simeq  \mathcal{L}^{n}_{2,a,ch} $.
\epf 

Now if $n=3$, we define two families of representations of $K$ on the vector space $V$ determined by the following action, for all $a \in \ku^{\times}$:
\begin{align*}
&\mathcal{L}^3_{1,a}: & & g \mapsto \begin{pmatrix}
\jordan_{2}(a) & 0  \\
0 & -a \\
\end{pmatrix}, & &
x_1 \mapsto \begin{pmatrix} 
0 & 0\\
e_2^t & 0  \end{pmatrix}, & &
x_2 \mapsto \begin{pmatrix} 
0 & 0 \\
ae_1^t & 0  \end{pmatrix}\\
&\mathcal{L}^3_{2,a}: & &  g \mapsto \begin{pmatrix}
\jordan_{2}(a) & 0  \\
0 & -a   \\ \end{pmatrix}, & & 
x_1 \mapsto \begin{pmatrix} 
0 & e_1\\
0 & 0 \\ \end{pmatrix}, & &
x_2 \mapsto \begin{pmatrix} 
0 & -ae_2\\
0 & 0 \end{pmatrix}. 
 \end{align*}

Clearly these are pairwise non-isomorphic indecomposable $K$-modules.
By a similar argument to the Lemma \ref{lema:two-blocks-p=1-super-result} we can prove the following result.

\begin{lema} If $V$ is indecomposable, then either $V \simeq\mathcal{L}^3_{1,a}$ or $V \simeq\mathcal{L}^3_{2,a}$. \qed
\end{lema}

\begin{obs}\label{dualL} Clearly, $({\mathcal{L}^3_{2,a}})^{*}\simeq \mathcal{L}^3_{1,a^{-1}}$, 
$\mathcal{L}^3_{2,1}\otimes \Bbbk_a\simeq \mathcal{L}^3_{2,a}\simeq \Bbbk_a\otimes \mathcal{L}^3_{2,1}$ and  $\mathcal{L}^3_{1,1}\otimes \Bbbk_{-a} \simeq \mathcal{L}^3_{1,a}\simeq \Bbbk_{-a}\otimes \mathcal{L}^3_{1,1}$, for all $a\in \Bbbk^{\times}$.
\end{obs}

\subsubsection{Tensor product of 3-dimensional indecomposable $K$-modules} \label{subsub-tensor}
Let $u_{ij}:=e_i\otimes e_j$, $i,j\in \I_	3$. Consider the basis $\base=\{v_i\,:i\in \I_9\}$ of $\mathcal{L}^3_{2,1}\otimes \mathcal{L}^3_{2,1}$, where
\begin{align*}
	&v_1= u_{31},\quad v_2=-u_{32},& &v_3=u_{13},& & v_4=-u_{23},& &v_5=u_{33},&\\
	&v_6=\frac{1}{2}(u_{12}-u_{21}),& &v_7=2u_{11},& & v_8= u_{12}+u_{21},& &v_9=-u_{12}+u_{22},&
\end{align*}
which give us the normal Jordan form of $g$. In the basis $\base$, the actions of the $g$, $x_1$ and $x_2$ are determined by \eqref{rep:general} where 
\begin{align*}
&A =\begin{pmatrix}
	\jordan_{2}(-1) & 0 \\
	0 & \jordan_{2}(-1)   
\end{pmatrix}, & 
&B = \begin{pmatrix}
	\id_{2} & 0  \\
	0 & \jordan_{3}(1)  
\end{pmatrix},&\\[.3em]
&C =\begin{pmatrix}
-e_1+e_3 & 0 & 0 & 0 &0 \\
\end{pmatrix}, & 
&E =\begin{pmatrix}
-e_2+e_4 & 0 & 0 & 0 &0 \\
\end{pmatrix}, & \\[.3em]
&D =\frac{1}{2}\begin{pmatrix}
0  \\
-2(e_2+e_4)^t \\
(e_1+e_3-e_4)^t\\ 
-(e_2-e_4)^t\\ 
0 \\  \end{pmatrix}, & 
&F =\begin{pmatrix}
0 \\
(e_1+e_2-e_3+2e_4)^t \\
0\\
-\frac{1}{2}(e_1-e_2+e_3-2e_4)^t \\
(e_2+e_4)^t \\ \end{pmatrix}. & 
\end{align*}

\begin{prop}\label{prop:tensorinde} 
	$T=\mathcal{L}^3_{2,1} \otimes \mathcal{L}^3_{2,1}$ is an indecomposable $K$-module.
\end{prop}
\pf Suppose that $T=U\oplus W$, with $U$ and $W$  non-trivial submodules of $T$. Note that $\ker x^2_2=\ker x_{21}=\ker x_2x_1=T-\la v_5\ra$, 
$ x_{21}T=\la v_7\ra$, $ x_2x_1T=\la v_6\ra$ and $ x^2_2T=\la v_6,v_8\ra$. We can assume $x_{21}U=\{0\}$ and $ x_{21}W=\la v_7\ra$. Let $u\in U$ and $w\in W$ such that $v_5=u+w$.
From $\la x_{21}w\ra =\la x_{21}v_5\ra=\la v_7\ra$ follows that $v_7\in W$. Similarly, from $ x_2x_1T=\la v_6\ra$ follows $v_6\in W$ and from $ x^2_2T=\la v_6,v_8\ra$ follows $v_8\in W$. Thus,  $W'=\la v_6,v_7,v_8 \ra \subset W$. Since $x_1u\in W'$, we obtain $u\in \ker x_1=\la v_1-v_3,v_6,v_7,v_8,v_9 \ra $. Hence, $x_2u\in \la v_6\ra$ and whence $u\in \ker x_2=\la v_6,v_7,v_8,v_9 \ra$. Since $u\notin W'$, it follows that $(g-\id)u=\alpha v_7+\beta v_8$ with $\alpha,\beta\in \ku$ and $\beta\neq 0$. Then $0\neq (g-\id)u\in U\cap W$ which is a contradiction.\epf

\begin{cor}\label{cor:tensor-geral} Let $a,b\in \Bbbk^{\times}$. Then
	\begin{enumerate}[leftmargin=*,label=\rm{(\roman*)}]
\item \label{item:cor:412-1} $ \mathcal{L}^{3}_{2,a}\otimes \mathcal{L}^{3}_{2,b}$ is an indecomposable $K$-module.
\item \label{item:cor:412-2} $\mathcal{L}^{3}_{1,a}\otimes \mathcal{L}^{3}_{1,b}$ is an indecomposable $K$-module.
	\end{enumerate}
\end{cor}
 
\pf By Remark \ref{dualL}, $\mathcal{L}^{3}_{2,a}\otimes \mathcal{L}^{3}_{2,b} \simeq \Bbbk_{ab}\otimes \mathcal{L}^{3}_{2,1}\otimes \mathcal{L}^{3}_{2,1}$. Hence \ref{item:cor:412-1} follows from Proposition \ref{prop:tensorinde}. By Remark \ref{dualL},  $\mathcal{L}^{3}_{1,a}\otimes \mathcal{L}^{3}_{1,b}\simeq (\mathcal{L}^3_{2,a^{-1}})^{*}\otimes (\mathcal{L}^3_{2,b^{-1}})^{*}\simeq (\mathcal{L}^3_{2,b^{-1}}\otimes \mathcal{L}^3_{2,a^{-1}})^{*}$. Thus, \ref{item:cor:412-2} follows from \ref{item:cor:412-1}.
\epf

Fix a basis $\{e_1,e_2,e_3\}$ of $\mathcal{L}^{3}_{2,1}$, a basis $\{\tilde{e}_1,\tilde{e}_2,\tilde{e}_3\}$ of  $\mathcal{L}^{3}_{1,1}$ and $u_{ij}:=e_i\otimes \tilde{e}_j$, $i,j\in \I_3$.  As above, the basis $\base=\{v_i\,: i\in \I_9\}$ of $\mathcal{L}^3_{2,1}\otimes \mathcal{L}^3_{1,1}$ with
\begin{align*}
&v_1= u_{31},\quad v_2=-u_{32},& &v_3=u_{13},& & v_4=-u_{23},& &v_5=u_{33},&\\
&v_6=\frac{1}{2}(u_{12}-u_{21}),& &v_7=2u_{11},& & v_8= u_{12}+u_{21},& &v_9=-u_{12}+u_{22},&
\end{align*}
give us the normal Jordan form of $g$ and the actions of the $g$, $x_1$ and $x_2$ are determined by \eqref{rep:general} where 
\begin{align*}
&A =\begin{pmatrix}
\jordan_{2}(-1) & 0 \\
0 & \jordan_{2}(-1)   
\end{pmatrix}, & 
&B = \begin{pmatrix}
\id_{2} & 0  \\
0 & \jordan_{3}(1)  
\end{pmatrix},&\\[.3em]
&C =\begin{pmatrix}
e_3 &\frac{1}{2} e_3  & 0 & e_3 & -e_4 \\
\end{pmatrix}, & 
&E =\begin{pmatrix}
e_4 & -\frac{1}{2}(e_3-e_4)  & 2e_3 & e_3-e_4 & 0 \\
\end{pmatrix}, & \\[.3em]
&D =\frac{1}{2}\begin{pmatrix}
2e_2^t \\
-2e_2^t \\
e_1^t  \\
-e_2^t\\
0\\
\end{pmatrix}, & 
&F =\begin{pmatrix}
-e_1^t\\
(e_1+e_2)^t  \\
0  \\
-\frac{1}{2}(e_1-e_2)^t  \\
e_2^t \\
\end{pmatrix}. & 
\end{align*}

\begin{prop}\label{prop:tensorcase2} 
	$\mathcal{L}^3_{2,1}\otimes \mathcal{L}^3_{1,1}\simeq U\oplus W$, where $U$ is an indecomposable $K$-module of dimension $8$ and $W=\ku_{1}$.
\end{prop}
\pf Consider $\mathfrak{C}=\{v_1, v_2, v_3,v_4,v_5-v_6,v_7,v_8,v_6+v_9\}$ which is a linearly independent set. 
Notice that $U=\la \mathfrak{C}\ra$ and $W=\left\langle v_5-2v_6-\frac{1}{2}v_7 \right\rangle$ are $K$-submodules of $\mathcal{L}^3_{2,1}\otimes \mathcal{L}^3_{1,1}$. Moreover  $\mathcal{L}^3_{2,1}\otimes \mathcal{L}^3_{1,1}\simeq U\oplus W$ and $W=\ku_{1}$.

Let $U_1$ and $U_2$ be non-trivial $K$-submodules of $U$ with  $U=U_1\oplus U_2$. Suppose that $v_3\in U_1$ and let $u\in U_2$ a vector with coordinates $(\lambda_i)$, $i\in \I_8$, in the basis $\mathfrak{C}$. From $x_2^2u=-\lambda_1v_3\in U_1\cap U_2$ we see that $\lambda_1=0$. Similarly, applying $x^2_2g,x_1$ and $x_2$ on $u$ we conclude that $\lambda_i=0$, for all $i\in \I_8$. Then $U_2=0$,  a contradiction. 
Hence, $v_3= u_1+u_2$ with $u_i\in U_i$ and $u_i\neq 0$ for $i\in\I_2$. Thus $u_1,u_2\in V^{-1}_g=\langle v_1,v_3\rangle$. Since $v_3\notin U_1$, it follows that $x_2^2\cdot u_1\neq 0$. Then $x_2^2\cdot u_1=-x_2^2\cdot u_2 \in \la v_3\ra$ and $v_3\in U_1\cap U_2$,   a contradiction.\epf

\begin{obs}
Since $\dim U^{\pm 1}=2$, $U\not\simeq \mathcal{L}^8_{i,a,b}$, for any $a\in \ku^{\times}$, $b\in \ku$, $i\in \I_2$. 
\end{obs}

\end{document}